%% file: ex_article.tex
\documentclass[hidelinks,onefignum,onetabnum]{siamart220329}

\input{ex_shared}

\ifpdf
\hypersetup{
  pdftitle={An Example Article},
  pdfauthor={D. Doe, P. T. Frank, and J. E. Smith}
}
\fi

\def \ihat{\mathbf{\hat{i}}}
\def \jhat{\mathbf{\hat{j}}}
\def \khat{\mathbf{\hat{k}}}
\def \ahat{\hat{\alpha}}
\def \bhat{\hat{\beta}}
\def \v{\mathbf{v}}
\def \r{\mathbf{r}}

\def \m{\mathbf{m}}
\def \khat{\mathbf{\hat{k}}}

\def \rp{\mathbf{r}_{p}}
\def \rq{\mathbf{r}_{q}}

\def \nmp{\mathbf{n}}

\def \np{\mathbf{n}_p}

\def \n{\mathbf{n}}

\def \K{K}
\def \V{V}
\def \Kp{K^{'}}
\def \D{D}

\newcommand{\R}{{\mathbb{R}}}

\begin{document}
\raggedbottom
\maketitle
\begin{abstract}
A method for the analytical evaluation of layer potentials arising in the collocation boundary element method for the Laplace and Helmholtz equation is developed for piecewise flat boundary elements with polynomial shape functions. 
The method is based on dimension-reduction via the divergence theorem and a Recursive scheme for evaluating the resulting line Integrals for Polynomial Elements (RIPE). 
It is used to evaluate single layer, double layer, adjoint double layer, and hypersingular potentials, for both the Laplace and the Helmholtz kernels. 
It naturally supports nearly singular, singular, and hypersingular integrals under a single framework without separate modifications.
The developed framework exhibits accuracy and efficiency. 
\end{abstract}

\begin{keywords}
Boundary element method,  singular integrals, analytical quadrature
\end{keywords}

\begin{MSCcodes}
65R20, 65N35, 65N38
\end{MSCcodes}

\section{Introduction}
Accurate, efficient, and error-controlled numerical evaluation of layer potentials is needed to build the linear systems that must be solved in boundary element methods (BEM). 
This is nontrivial because the integrands are not simple analytic functions but are singular or nearly singular. 
Quadrature schemes which are effective for integrating band-limited polynomials, e.g. Gauss-Legendre quadrature, are known to produce inaccurate results when the evaluation point is close to the element. 
Many techniques have been developed over the years to accurately evaluate boundary integrals in such cases~\cite{hayami1988quadrature, guiggiani1992general, hackbusch1994numerical, johnston2007sinh, lenoir2012evaluation, montanelli2021computing, zhu2022high, adelman2016computation, klockner2013quadrature, wala2019fast, wala2020optimization}. 
Ref. \cite{montanelli2021computing} provides a recent extensive survey on this subject. 
The approaches developed include 
singularity cancellation using coordinate transforms~\cite{hayami1988quadrature, hackbusch1994numerical, johnston2007sinh}, 
singularity subtraction~\cite{guiggiani1992general}, 
dimension reduction~\cite{lenoir2012evaluation, zhu2022high}, 
adaptive subdivision~\cite{adelman2016computation}, and quadrature by expansion~\cite{klockner2013quadrature,wala2019fast,wala2020optimization}. 
The PART method~\cite{hayami1988quadrature} applies coordinate transforms which reduce the effect of the diverging integrand via the Jacobian of the transform. 
\cite{johnston2007sinh} introduced another coordinate transform approach, using $\sinh$ transforms. 
Singularity subtraction methods~\cite{guiggiani1992general} split the integrand into a singular part and a regular part, where the former is evaluated analytically and the latter via Gauss-Legendre quadrature. 
Another approach for evaluating the layer potentials is to derive analytical expressions of the integrals. 
Newman developed a method to evaluate layer potentials for the Laplace kernel on quadrilateral elements for shape functions of arbitrary order~\cite{newman1986distributions}. 
Lenoir and Salles developed a semi-analytical method using dimension-reduction via the divergence theorem for constant and linear elements in Galerkin BEM~\cite{lenoir2012evaluation, lenoir2012exact}, where the layer potential associated with the singular part of the Helmholtz kernel is evaluated analytically, but the regular part is evaluated via conventional quadrature. 
Recently, Zhu and Veerapaneni~\cite{zhu2022high} introduced a method for Laplace layer potentials on high-order curved elements using dimension reduction via Stokes' theorem and quaternion algebra. 

We propose a method to analytically evaluate all four layer potentials arising in collocation BEM for the Helmholtz and Laplace kernels over piecewise flat boundary elements with high-order polynomial shape functions. 
The method is based on reduction of the surface integrals to line integrals via the divergence theorem in the element plane and efficient computation of the resulting line integrals using recursions which are derived via the introduction of auxiliary vector fields. 
In comparison to \cite{lenoir2012evaluation} we treat both the Laplace and Helmholtz kernels within a single framework, and are able to analytically evaluate integrals with general order polynomial shape functions via recursions. 
The method in \cite{zhu2022high} is most related to ours in that it is based on dimension reduction and supports polynomial elements. 
However our work differs from \cite{zhu2022high} in a number of ways: 
(1) it is analytical and solely uses recursions, in contrast to the use of Gauss-Legendre quadrature for the line integrals, 
(2) it is derived for both the Laplace and Helmholtz kernels and for all four basic layer potentials, i.e. single layer, double layer, adjoint double layer, and hypersingular potential, 
(3) it is restricted to flat elements by relying on the divergence theorem, whereas \cite{zhu2022high} is based on Stokes' theorem for manifolds and can handle curved elements, and 
(4) where \cite{zhu2022high} involves an $O(p^6)$ coordinate transform step,  our method has an overall complexity of $O(p^3)$ for the Laplace case, with $p$ the polynomial order. 
We remark that the divergence theorem in the flat element plane can be derived from Stokes' theorem, so the method could have been equivalently presented in terms of Stokes' theorem. 
 
One of the benefits of our recursion-based approach is error control, as some of the recursions may be truncated as soon as a prescribed threshold is achieved, a feature that is missing in fixed-order Gauss-Legendre quadrature. 
Error control is important when quadrature is part of a more complicated solver with multiple components, e.g. the fast multipole method (FMM)~\cite{greengard1987fast} accelerated BEM solved via an iterative method, where the numerical error arising from different parts should all be kept consistent to avoid redundant computation. While adaptive error control can  be achieved by classical methods e.g. Gauss-Kronrod quadrature, such methods are not as efficient as our method, as seen in \S 5. 
The simplest case of our method for piecewise constant elements was initially developed by Gumerov and Duraiswami~\cite{gumerov2021analytical} and used in a production collocation BEM solver for the Helmholtz equation~\cite{gumerov2021fast}.

\section{Boundary element method and layer potentials}
The boundary element method is extensively used for numerical solution of partial differential equations, e.g. the Helmholtz equation and the Laplace equation, respectively given by
\begin{equation}\begin{aligned}
 - k^2 u(\r) - \nabla^2 u(\r) = f(\r), \quad  - \nabla^2 u(\r) = f(\r),
\end{aligned}\label{eq:HelmholtzEqn} \end{equation}
with wavenumber $k$, field $u$, and source $f$. 
The weak form of \cref{eq:HelmholtzEqn} can be written in terms of single- and double layer potentials $V$, $K$~\cite{sauter2010boundary}: 
\begin{equation}\begin{aligned}
\{( c_p \gamma_{0,p} + \K\gamma_{0,q} - \V\gamma_{1,q}) u\}(\rp) &= \{N_0f\}(\rp),
\end{aligned}\label{eq:BIE_BF_withRobinBC}
\end{equation}
with $c_p=1/2$, $\gamma_{0}$ and $\gamma_{1}$ the boundary trace and normal derivative operators, $N_0$ the Newton potential operator, and layer potentials defined over the boundary $\Gamma$ as:
\begin{equation}\begin{aligned}
\{\V\psi\}(\rp) \equiv \int_{\rq \in \Gamma} G(\mathbf{r}_p,\mathbf{r}_q) \psi(\rq) d\Gamma, \quad
\{ \K \phi\}(\rp) \equiv \int_{\rq \in \Gamma} \frac{\partial G(\mathbf{r}_p,\mathbf{r}_q)}{\partial \nmp_q} \phi(\rq) d\Gamma,\\
\end{aligned} \end{equation}
where $G(\rp,\rq)$ is the Helmholtz or Laplace Green function:
\begin{equation}\begin{aligned}
G_\mathrm{H}(\rp,\rq) = \frac{e^{ikr}}{4\pi r}, \quad G_\mathrm{L}(\rp,\rq) = \frac{1}{4\pi r}, \quad r \equiv |\mathbf{r}_q-\mathbf{r}_p|.
\end{aligned} \end{equation}
Some BEM formulations~\cite{burton1971application, gumerov2013fast} use the normal derivative form of the boundary integral equation with the adjoint double layer potential $\Kp$ and hypersingular potential $D$.
 
\begin{equation}\begin{aligned}
\!\!\{  \Kp \psi\}(\rp) \equiv \int_{\rq \in \Gamma} \!\!\!\!\frac{\partial G(\mathbf{r}_p,\mathbf{r}_q)}{\partial \nmp_p} \psi(\rq) d\Gamma, \!\!\quad
\{  \D \phi\}(\rp) \equiv \int_{\rq \in \Gamma}\!\!\!\! \frac{-\partial^2 G(\mathbf{r}_p,\mathbf{r}_q)}{\partial \nmp_p \partial \nmp_q} \phi(\rq) d\Gamma. \\
\end{aligned} \end{equation}
In the collocation BEM the boundary $\Gamma$ is discretized into elements, typically triangular, and the layer potential integrals over these elements are evaluated. 
The densities $\psi$, $\phi$ are approximated via local polynomial functions (also called \emph{shape functions}). 

The integrals involved in evaluation of the layer potentials over elements can be classified depending on the distance of the evaluation point from the element and the degree of singularity as: 
(1) {\em far-field} where the evaluation points are far from the element, and typically clustered via the FMM for acceleration, 
(2) {\em intermediate} distance where the FMM is not used but standard Gauss-Legendre quadrature is still effective, 
(3) {\em near-field} where Gauss-Legendre quadrature becomes inaccurate and special treatment for nearly singular integrals is needed, and 
(4) evaluation points on the elements where the integrand is either {\em singular or hypersingular} and other special techniques are needed. 
Techniques have been developed for each of these categories. 
For the far-field, the authors have developed another recursion-based method for the efficient evaluation of integrals of multipoles arising in the FMM~\cite{gumerov2021fast}. 
The quadrature by expansion method has also been used with the FMM~\cite{wala2019fast,wala2020optimization}. 

\section{Problem statement}
\begin{figure}[htbp]
    \centering
    \includegraphics[width=10cm]{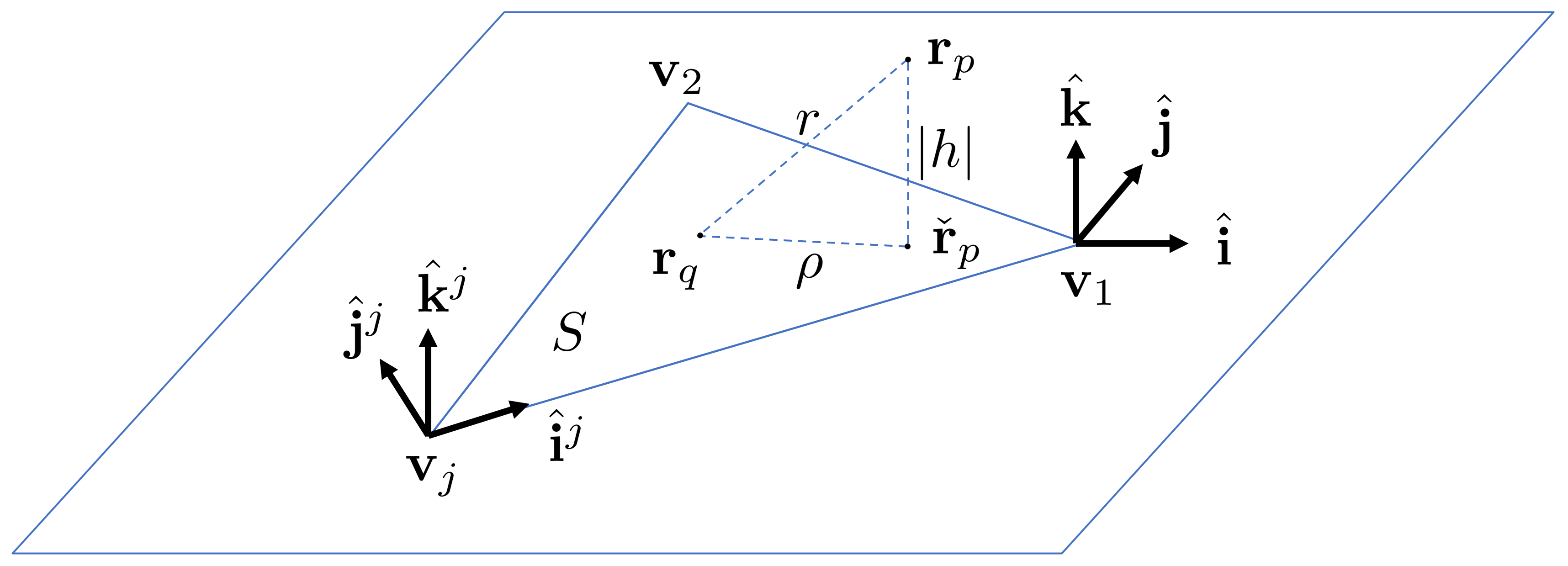}
    \caption{A boundary element and its associated coordinate frames.
     We have the \emph{element coordinate frame} $(\ihat,\jhat,\khat)$ and \emph{edge coordinate frame} $(\ihat^j, \jhat^j, \khat^j)$ for each edge with index $j$. The origin of the element frame is set to $\v_1$.
}
    \label{fig:elementAndCoordinateFrames}
    \end{figure}

Assume a flat triangle element $S$ with vertices $\v_1$, $\v_2$, and $\v_3$ in $\R^3$, where the \emph{element coordinate frame} $(\ihat,\jhat,\khat)$ has its origin in the element plane, $\ihat$ and $\jhat$ are unit vectors lying in the element plane, and $\khat = \ihat \times \jhat = \n_{\mathrm{face}}$, where $\n_{\mathrm{face}}$ is the unit normal to the element. 
We also define the \emph{edge coordinate frame} $(\ihat^j, \jhat^j, \khat^j)$ for edge $j$ as: 
$\ihat^j = (\v_{(j+1)\%3}-\v_{j})/|\v_{(j+1)\%3}-\v_{j}|$, 
$\khat^j = \n_{\mathrm{face}}$, and
$\jhat^j = \khat^j \times \ihat^j \equiv -\n^j$.
$\n^1$, $\n^2$, and $\n^3$ are the unit outward edge normal vectors in the element plane for edges $(\v_1,\v_2)$, $(\v_2,\v_3)$, and $(\v_3,\v_1)$, respectively.
The coordinate frames associated with the boundary element are shown in \cref{fig:elementAndCoordinateFrames}. In the following, the superscript $j$ refers to the edge index and is not an exponent. 
$\rp$ and $\check{\r}_p$ are respectively defined as the observation point and its orthogonal projection on the element plane. 
The points $\rq$ and $\rp$ are represented as $\rq = (x,y,z)^T$ and $\rp = (x_p,y_p,z_p)^T$ in the element frame. 
Also define $r \equiv |\mathbf{r}_q-\mathbf{r}_p|$, $\boldsymbol{\rho} \equiv \r_q-\check{\mathbf{r}}_p $, $\rho \equiv |\boldsymbol{\rho} |$, $h \equiv \khat \cdot (\rp-\check{\mathbf{r}}_p )$, $x_d \equiv x - x_{p}$, $y_d \equiv y - y_{p}$, and $z_d \equiv z - z_{p}$.

In this setup, a point in space can be expressed as:
\begin{equation}
\r = \v_1 + x \ihat + y\jhat + z\khat = \v_1 + u\left(\v_2-\v_1\right)+v\left(\v_3-\v_1\right) + w \n_{\mathrm{face}},
\end{equation}
where $(u,v)$ are the coordinates in the two-dimensional reference frame
with $0 \le u,v \le 1$ and $u+v \le 1$, and $(x,y,z)$ are the coordinates in the element frame. 
Projection from the reference frame to the element frame is a linear transform: 
\begin{equation}
(x,y)^T = \Xi (u,v)^T,
\label{eq:xy_and_uv}\end{equation}
with $\Xi$ the $(2\times2)$ transformation matrix. 
Points in the element frame are expressed by their coordinates in the edge frame as: 
\begin{equation}
(x,y,z)^T
=  B (x^j,y^j,z^j)^T + \mathbf{b},\quad 
B \equiv \left(\begin{array}{ccc}
\alpha_x & \beta_x & 0 \\
\alpha_y & \beta_y & 0 \\
0 & 0 & 1
\end{array}\right), \quad \mathbf{b} \equiv \left(\begin{array}{l}
\gamma_x \\
\gamma_y \\
0
\end{array}\right).
\end{equation}
Identifying the entries of $\Xi$, $B$, and $\mathbf{b}$ for a given element is straightforward. 

Shape functions $N(u,v)$ for $p_s$-th order elements have the following form:
\begin{equation}
\begin{aligned}
N_{p_s}(u,v) & = \sum_{b'=0}^{p_s} \sum_{c'=0}^{p_s-b'} a_{b',c'} u^{b'} v^{c'} = \sum_{b=0}^{p_s} \sum_{c=0}^{p_s-b} A_{b,c} x^b y^c.
\end{aligned}\label{eq:shapeFncs}
\end{equation}
Shape functions are specified by the $a_{b',c'}$ coefficients in the $(u,v)$ reference frame. 
The coefficients $A_{b,c}$ for the element frame expression can be obtained by substituting \cref{eq:xy_and_uv} into \cref{eq:shapeFncs} and solving a linear system.
Since the coefficients $a_{b',c'}$ and $A_{b,c}$ only interact if $b+c=b'+c'$ due to \cref{eq:xy_and_uv}, this can be achieved by solving $p_s+1$ small linear systems, 
and each of these can be solved via efficient stabilized algorithms for Vandermonde systems~\cite{bjorck1970solution}, resulting in a complexity of $O(p_s^3)$ for this stage.

\section{Method: Recursive Integrals for Polynomial Elements (RIPE)}
\subsection{Dimension reduction via divergence theorem} 
\begin{figure*}[htbp]
    \centering
    \includegraphics[width=10cm]{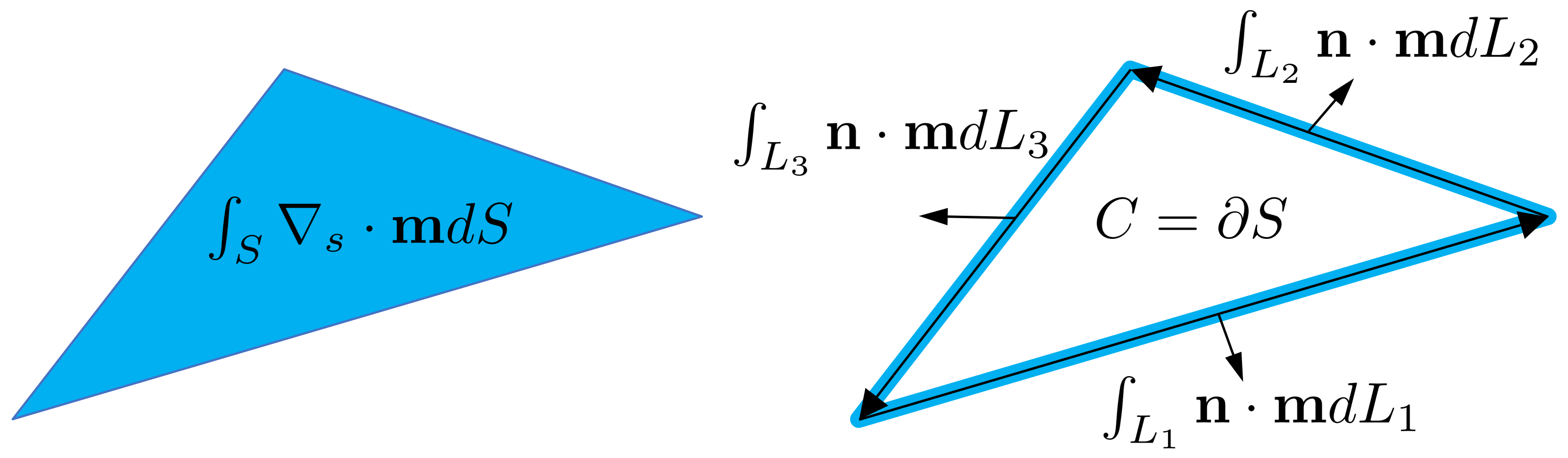}
    \caption{Left: the surface integral over element $S$ of a function which can be written as the surface divergence of a vector field $\m$. Right: after the application of the divergence theorem, the surface integral over $S$ is converted into a contour integral over the contour $C=\partial S$ and further into line integrals due to the straightness of the edges.}
    \label{fig:S_to_dS_via_Gauss_basic}
    \end{figure*}

\begin{figure*}[htbp]
    \centering
    \includegraphics[width=10cm]{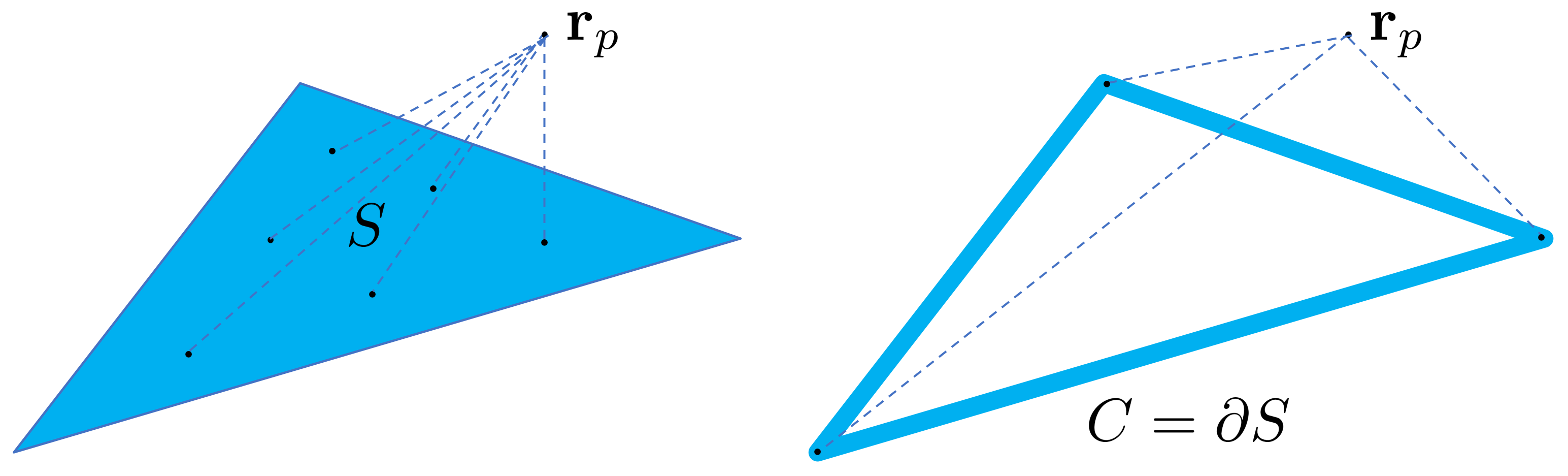}
    \caption{Conventional Gauss-Legendre quadrature requires  evaluation of the integrand at multiple nodes on the element $S$, which may be near or coincide with the evaluation point, or collocation node, $\rp$. After application of the divergence theorem, the surface integral over $S$ is converted into an integral over the contour $C=\partial S$ and further into line integrals, which in the end results in the evaluation of analytical expressions at the endpoints of the edges, i.e. at the vertices.}
    \label{fig:S_to_dS_via_Gauss}
    \end{figure*}

We assume flat triangular elements and reduce the surface integral to a contour integral via the divergence theorem, and discretize the contour integral into a sum of line integrals over the straight edges. 
We seek vector fields $\mathbf{m}_i$ whose surface divergence in the element plane gives the integrand in question, i.e. a polynomial shape function multiplied with the Green function: 
 \begin{equation}\begin{aligned}
\sum_i A_i \nabla_s \cdot \m_i = N(u(x,y),v(x,y)) G(\rp,\rq),
 \end{aligned}\end{equation}
then we can apply the divergence theorem for the evaluation of the element-wise layer potentials of interest.
This is depicted in \cref{fig:S_to_dS_via_Gauss_basic}.
The single layer potential $\V_e$, double layer potential $\K_e$, adjoint double layer potential $\Kp_e$, and hypersingular potential $\D_e$ can be all then written as:
\begin{equation}\label{eq:layerPotentials}\begin{aligned}
\!\!\!\!\V_e &= \!\!\int_{S_q} \!\!G(\rp,\rq) N(\rq) dS_q
= \!\!\int_{S_q} \!\!\sum_{i} \!\!A_i \nabla_s \cdot \m_i dS_q = \sum_{i} A_i  \int_{S_q} \nabla_s \cdot \m_i dS_q\\
&= \sum_{i} A_i \left( \oint_{C} \n_{C} \cdot \m_i dC  + T_{\mathrm{sing}}^{(i)} \right) 
= \sum_{i}  A_i \left( \sum_{j}  \int_{L_j} \n_{j} \cdot \m_i dL_j + T_{\mathrm{sing}}^{(i)} \right),\\
\K_e &= \int_{S_q} \frac{\partial G(\rp,\rq)}{\partial \nmp_q} N(\rq) dS_q 
=  \sum_{i} A_i   \left( \sum_{j} \frac{\partial}{\partial z} \int_{L_j} \n_{j} \cdot \m_i dL_j  + \frac{\partial T_{\mathrm{sing}}^{(i)} }{\partial z} \right),\\
\Kp_e &= \int_{S_q} \frac{\partial G(\rp,\rq)}{\partial \nmp_p} N(\rq) dS_q 
= \sum_{i} A_i  \left( \sum_{j}  \frac{\partial}{\partial \np}  \int_{L_j} \n_{j} \cdot \m_i dL_j  +\frac{\partial T_{\mathrm{sing}}^{(i)}}{\partial \np}   \right),\\
-\D_e &= \int_{S_q} \frac{\partial^2 G(\rp,\rq)}{\partial \nmp_p\partial \nmp_q} N(\rq) dS_q 
= \sum_{i} A_i  \left( \sum_{j}\frac{\partial^2}{\partial \np \partial z} \int_{L_j} \n_{j} \cdot \m_i dL_j  +\frac{\partial^2 T_{\mathrm{sing}}^{(i)} }{\partial \np \partial z}  \right),\\
\end{aligned}\end{equation}
with $\nabla_s=\nabla-\khat(\khat \cdot \nabla)$ the surface divergence operator in the element plane and $C=\partial S_q$ the contour of the element surface $S_q$. 
Note that we have to add the singularity term $T_{\mathrm{sing}}^{(i)}$ associated with $\m_i$ to the contour integral if  $\check{\mathbf{r}}_p$, the projection of $\rp$ onto the element plane, falls onto the element and if the integrand $\n_C \cdot \m_i$ has a singularity at $\rq=\check{\mathbf{r}}_p$ on the element. 
In this case, the singularity has to be removed when applying the divergence theorem, i.e. the contour integral over an infinitesimal circle $C_{\varepsilon}$ centered at the singularity on the element has to be subtracted from the contour integral over $S$: 
$T_{\mathrm{sing}}^{(i)} \equiv \lim_{\varepsilon \to 0} \oint_{C_{\varepsilon}} \n_{C_{\varepsilon}} \cdot \m_i dC_{\varepsilon}$.
Thus, the problem reduces to the evaluation of line integrals $\int_{L_j} \n_j \cdot \m_i dL_j$ associated with each edge and singularity terms associated with the element. 
One way to evaluate these line integrals is to rely on Gauss-Legendre quadrature, resulting in semi-analytical methods~\cite{zhu2022high}. 
If these line integrals can be computed analytically, they can be further reduced to the evaluation at the end-points of the edges due to the fundamental theorem of calculus. 
\cref{fig:S_to_dS_via_Gauss} depicts the contrast of the conventional approach of surface integrals based on Gauss-Legendre quadrature and the proposed RIPE method. 

To execute this strategy we need to find expressions for the vector fields $\m_i$, whose surface divergence equals the integrand $G(\rp,\rq) N(\rq)$, or, at least, the terms $\int_{L_j} \n_j \cdot \m_i dL_j$ and $T_{\mathrm{sing}}^{(i)}$. 
While the proposed method can be naturally generalized to polygonal elements, here we only consider triangular elements for their simplicity and efficiency. 
In the rest of the paper, ${\partial f}/{\partial x} = \partial_x f = f^{(x)}$ all denote partial derivatives of a function $f$ with respect to variable $x$.

\subsection{Laplace kernel with polynomials of arbitrary order}
\subsubsection{Single layer potential}
Our aim is to find vector fields whose surface divergence equals the integrand in question, i.e., a polynomial of $x$ and $y$ with degree up to $p_s$ multiplied with the Green function:
\begin{equation}\begin{aligned}
\nabla_s \cdot \mathbf{m}_{b,c} = \frac{ x^b y^c}{4 \pi r}, \quad b,c \ge 0, \quad b+c \le p_s. \\
\end{aligned}\label{eq:definition_of_m_bc_laplace}\end{equation}
We introduce auxiliary vector fields $\mathbf{m}_{b,c,a}$ which satisfy the following relation:
\begin{equation}\begin{aligned}
\nabla_s \cdot \mathbf{m}_{b,c,a} = \frac{ x^b y^c r^{a} }{4 \pi r}, \quad a,b,c \ge 0, \quad a+b+c \le p_s. \\
\end{aligned}\label{eq:definition_of_m_abc_laplace}\end{equation}
\begin{lemma}[Recurrence for auxiliary vector fields - Laplace kernel]\label{lemma:recurrence_b_c_laplace}
  If fields $\m_{b-1,c,a+1}$ and $\m_{b,c,a-1}$ satisfy \cref{eq:definition_of_m_abc_laplace}, then the field $\m_{b+1,c,a-1}$  given by:
   \begin{equation}\begin{aligned}
a \mathbf{m}_{b+1,c,a-1} = 
- b \mathbf{m}_{b-1,c,a+1}   + a x_p \mathbf{m}_{b,c,a-1} +  \left(\frac{x^{b} y^c r^a }{4 \pi}\ihat \right), \quad a\ge1, \\
\end{aligned}\label{eq:rec_m_bca_laplace}\end{equation}
 also satisfies \cref{eq:definition_of_m_abc_laplace} for the index tuple $(b+1,c,a-1)$.
\end{lemma}
\begin{proof}
 If vector fields $\m_{b-1,c,a+1}$ and $\m_{b,c,a-1}$ satisfy \cref{eq:definition_of_m_abc_laplace}, the surface divergence of the right hand side of \cref{eq:rec_m_bca_laplace} yields:
 \begin{equation}\begin{aligned}
&- b \nabla_s \cdot \mathbf{m}_{b-1,c,a+1}   + a x_p \nabla_s \cdot \mathbf{m}_{b,c,a-1} + \nabla_s \cdot \left(\frac{x^{b} y^c r^a }{4 \pi}\ihat \right) \\
&= -b \frac{x^{b-1}y^cr^a}{4\pi r } + a x_p \frac{x^{b}y^cr^{a-1}}{4\pi r } + y^{c}\frac{a x^{b+1}r^{a-1}- a x_p x^{b}  r^{a-1}  + b x^{b-1} r^{a+1}   }{4 \pi r} \\
&= a\frac{x^{b+1}y^cr^{a-1}}{4\pi r }, \\
\end{aligned}\end{equation}
hence, $\nabla_s \cdot \mathbf{m}_{b+1,c,a-1}= x^{b+1}y^cr^{a-1}/(4\pi r)$, which means that $\m_{b+1,c,a-1}$ satisfies \cref{eq:definition_of_m_abc_laplace} for the index tuple $(b+1,c,a-1)$.
\end{proof}
Let us define:
\begin{equation}\begin{aligned}
\xi_{b,c,a} \equiv \int \mathbf{m}_{b,c,a}  \cdot \n^j dx^j, \quad i_m \equiv \int r^m dx^j, \quad k_m \equiv -y_p^j \int \frac{r^m}{\rho^2} dx^j.
\end{aligned}\label{eq:definition_of_xi_abc_laplace}\end{equation}
From \cref{eq:rec_m_bca_laplace} it follows: 
\begin{equation}\begin{aligned}
  \xi_{b+1,c,a} = - \frac{b}{a+1} \xi_{b-1,c,a+2} + x_p \xi_{b,c,a} + \frac{\ihat \cdot \n^j }{(a+1)4\pi} \int x^{b} y^c r^{a+1} dx^j,\\
\end{aligned}\label{eq:rec_xi_bca_laplace_b}\end{equation}
and by symmetry, the following recurrence also holds:
\begin{equation}\begin{aligned}
\xi_{b,c+1,a} = - \frac{c}{a+1} \xi_{b,c-1,a+2} + y_p \xi_{b,c,a} + \frac{\jhat \cdot \n^j }{(a+1)4\pi} \int x^{b} y^c r^{a+1} dx^j.\\
\end{aligned}\label{eq:rec_xi_bca_laplace_c}\end{equation}
These recursions hold for $b\ge0$ and $c\ge0$. 
For $b=c=0$, the following field $\m_{0,0,a}$:
\begin{equation}\begin{aligned}
\m_{0,0,a} = \frac{r^{a+1} }{(a+1)4\pi\rho^2}\boldsymbol{\rho},
\end{aligned}\label{eq:m_0a0_laplace}\end{equation}
satisfies \cref{eq:definition_of_m_abc_laplace} and can be used to initiate the recursion. This yields:
\begin{equation}\begin{aligned}
\xi_{0,0,a} = \int \frac{r^{a+1} }{(a+1)4\pi\rho^2}\boldsymbol{\rho} \cdot \n^j dx^j =  \frac{y_p^j}{(a+1)4\pi} \int \frac{r^{a+1} }{\rho^2} dx^j = \frac{- k_{a+1}}{(a+1)4\pi}.
\end{aligned}\label{eq:rec_of_xi_00a_laplace}\end{equation}
 The integral $k_m$ and $i_m$ can be computed efficiently via the following recurrence: 
\begin{equation}\begin{aligned}
 k_{m+2} = (z_p^j)^2 k_m  - y_p^j i_m , \quad i_{m+2} = ((y_p^j)^2+(z_p^j)^2) \frac{m+2}{m+3} i_m + \frac{x_d^j r^{m+2}}{m+3},
\end{aligned}\label{eq:rec_k_i}\end{equation}
with closed form initial terms in \cref{appendix_kl_il}.
The following are used later:
\begin{equation}\begin{aligned}
\partial_{x_p^j} k_l = y_p^j \frac{r^l}{\rho^2},\quad
\partial_{y_p^j} k_l &= x_d^j \frac{r^l}{\rho^2} - l i_{l-2},\quad
\partial_{z_p^j} k_l = -l z_d^j k_{l-2},\\
\partial_{x_p^j} i_l = -r^l,\quad 
\partial_{y_p^j} i_l &= l y_p^j i_{l-2},\quad
\partial_{z_p^j} i_l = -l z_d^j i_{l-2}.\\
\end{aligned}\end{equation}
It remains to compute the last integral in \cref{eq:rec_xi_bca_laplace_b} and \cref{eq:rec_xi_bca_laplace_c} with form $\int x^b y^c r^a dx^j$.
Recall $x=\alpha_x x^j + \beta_x y^j + \gamma_x$ and $y=\alpha_y x^j + \beta_y y^j + \gamma_y$.
Also, $y^j=z^j=0$ for points on edge $j$. Following a coordinate transform from $x^j$ to $x^j_d$ this can be written as:
\begin{equation}\begin{aligned}
& \kappa_{b,c,a} \equiv \int x^{b}  y^c r^a   dx^j
 = \int (\alpha_x x^j_d-\hat{\alpha})^b   (\alpha_y x^j_d-\hat{\beta})^c r^a dx^j_d, \\
& \ahat \equiv -(\beta_x y^j+\gamma_x) - \alpha_x x_p^j, \quad 
\bhat \equiv - (\beta_y y^j+\gamma_y) - \alpha_y x_p^j, \quad \gamma^2 \equiv r^2-(x_d^j)^2.
 \end{aligned}\end{equation}
Later we use:
\begin{equation}\begin{aligned}
\ahat^{(x_p^j)} = - \alpha_x, \quad \bhat^{(x_p^j)} = -\alpha_y, \quad \ahat^{(y_p^j)} = \ahat^{(z_p^j)} = \bhat^{(y_p^j)} = \bhat^{(z_p^j)} = 0. \\
\end{aligned}\label{eq:derivatives_of_alphaHat_and_betaHat}\end{equation}

The function $\kappa_{b,c,a}$ satisfies the following recurrence relations:
\begin{equation}\label{eq:rec_of_kappa_bplus_cplus_laplace}\begin{aligned}
\kappa_{b+1, c+1, a}
&=\alpha_x\alpha_y \kappa_{b, c, a+2} -  \bhat \kappa_{b+1, c, a} -\ahat \kappa_{b,c+1,a} - \left(\ahat\bhat+\alpha_x\alpha_y\gamma^2\right) \kappa_{b, c, a},\\
\kappa_{b+2, c, a}
&=\alpha_x^2 \kappa_{b, c, a+2} - 2 \ahat \kappa_{b+1, c, a}-\left(\ahat^2+\alpha_x^2\gamma^2\right) \kappa_{b, c, a},\\
\kappa_{b, c+2, a}
&=\alpha_y^2 \kappa_{b, c, a+2} - 2 \bhat \kappa_{b, c+1, a}-\left(\bhat^2+\alpha_y^2\gamma^2\right) \kappa_{b, c, a}.\\
\end{aligned}\end{equation}

The following special cases are used to initiate the recursion:
\begin{equation}\begin{aligned}
&\kappa_{0,0,a} = i_{a},\quad
\kappa_{1, 0, a} = \alpha_x \int x^j_d r^a dx^j_d - \ahat i_a = \alpha_x \frac{r^{a+2}}{a+2} - \ahat \kappa_{0,0,a},\\
&\kappa_{0, 1, a} = \alpha_y \int x^j_d r^a dx^j_d - \bhat i_a = \alpha_y \frac{r^{a+2}}{a+2} - \bhat \kappa_{0,0,a}.\\
\end{aligned}\label{eq:rec_of_kappa_00a_01a_10a_11a_laplace}\end{equation}

Later, we will use the following relations:
\begin{equation}\begin{aligned}
\partial_{x_p^j} \kappa_{b,c,a} 
  &=  \partial_{x_p^j} \int x^{b} y^c r^{a} dx^j 
  =  \int x^{b} y^c \partial_{x_p^j} (r^{a}) dx^j 
  =  -\int x^{b} y^c \partial_{x^j} (r^{a}) dx^j \\
  &= - x^b y^c r^a + \alpha_x b \kappa_{b-1,c,a} + \alpha_y c \kappa_{b,c-1,a},\\
\partial_{y_p^j} \kappa_{b,c,a} &= - a y_d \kappa_{b,c,a-2},\quad \partial_{z} \kappa_{b,c,a} = a z_d \kappa_{b,c,a-2},\\
\partial^2_{z,x_p^j} \kappa_{b,c,a} 
  &= - x^b y^c a z_d r^{a-2} + \alpha_x b \kappa_{b-1,c,a}^{(z)} + \alpha_y c \kappa_{b,c-1,a}^{(z)}\\
  &= a z_d (- x^b y^c r^{a-2} + \alpha_x b \kappa_{b-1,c,a-2} + \alpha_y c \kappa_{b,c-1,a-2}).\\ 
\end{aligned}\end{equation}

\begin{figure}[htbp]
\centering
\includegraphics[width=4cm]{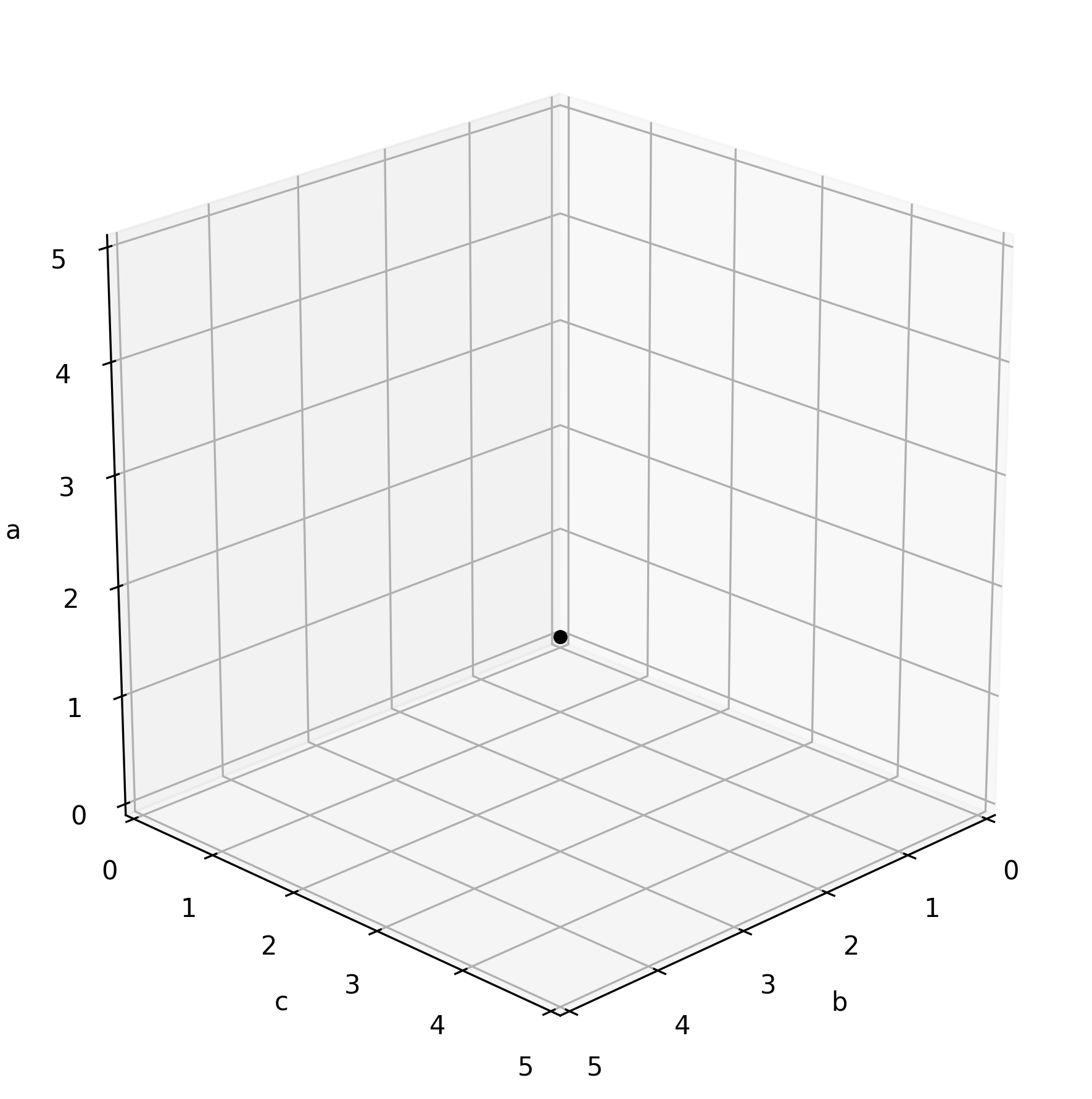}
\includegraphics[width=4cm]{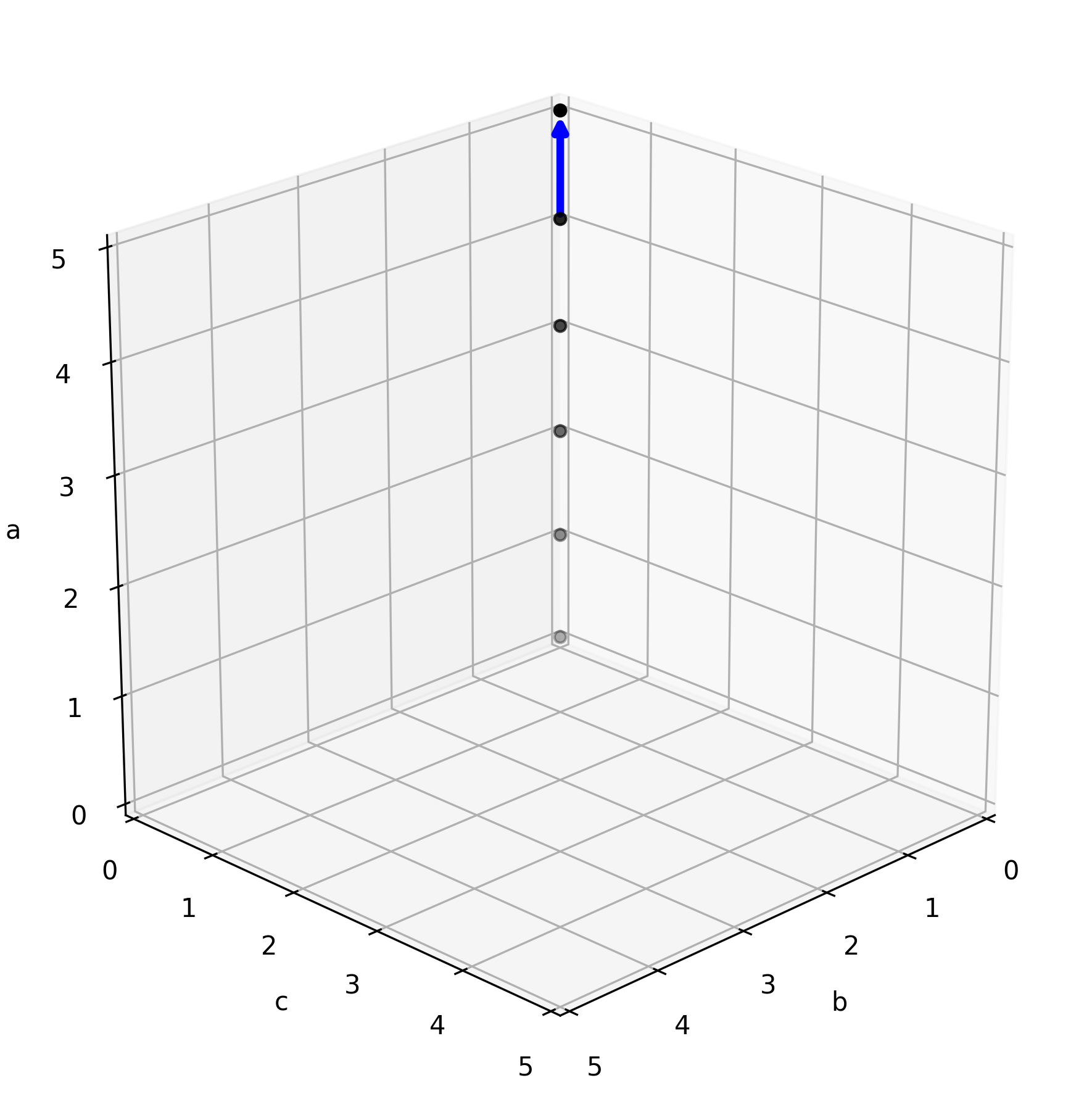}
\includegraphics[width=4cm]{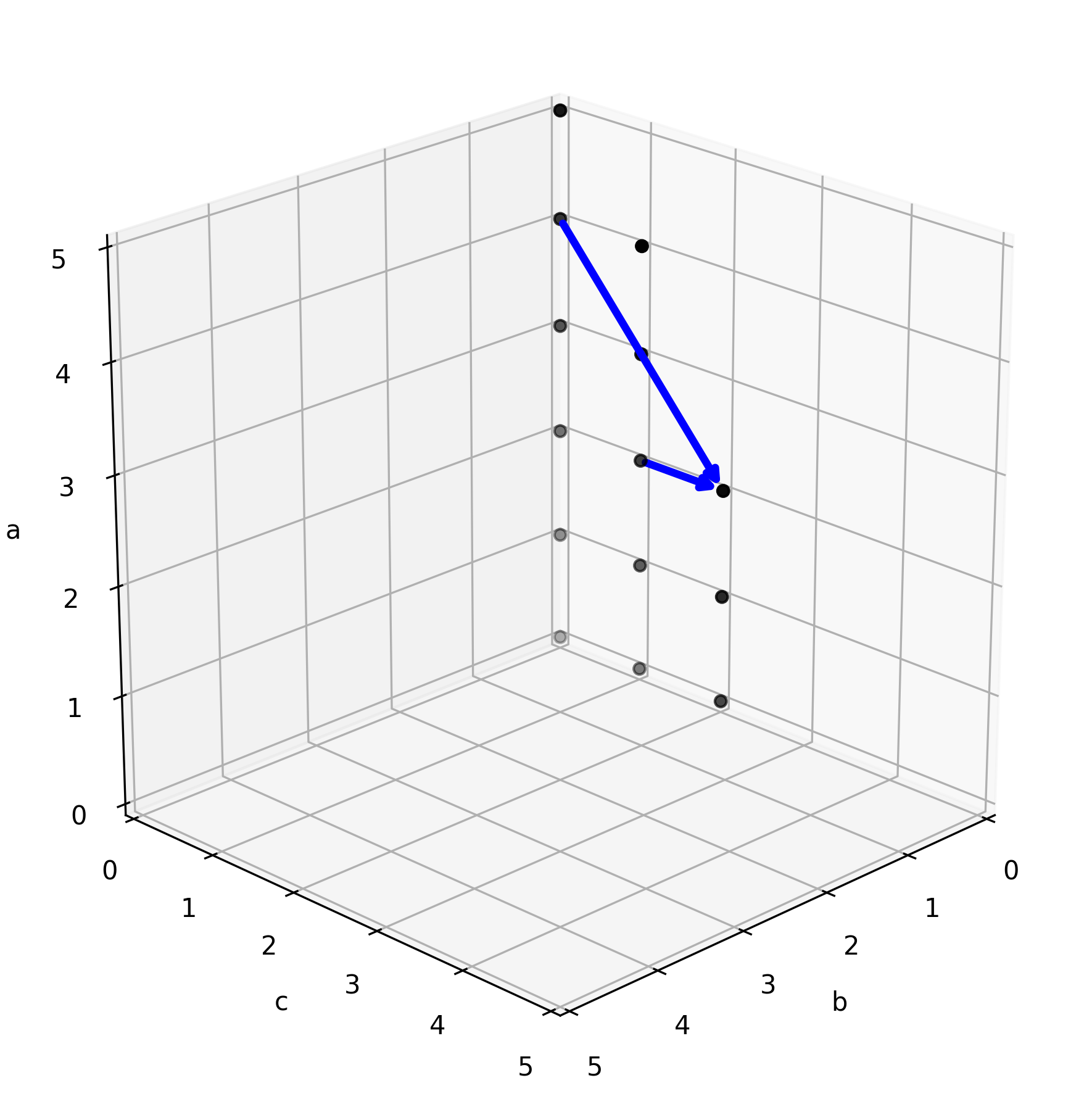}
\includegraphics[width=4cm]{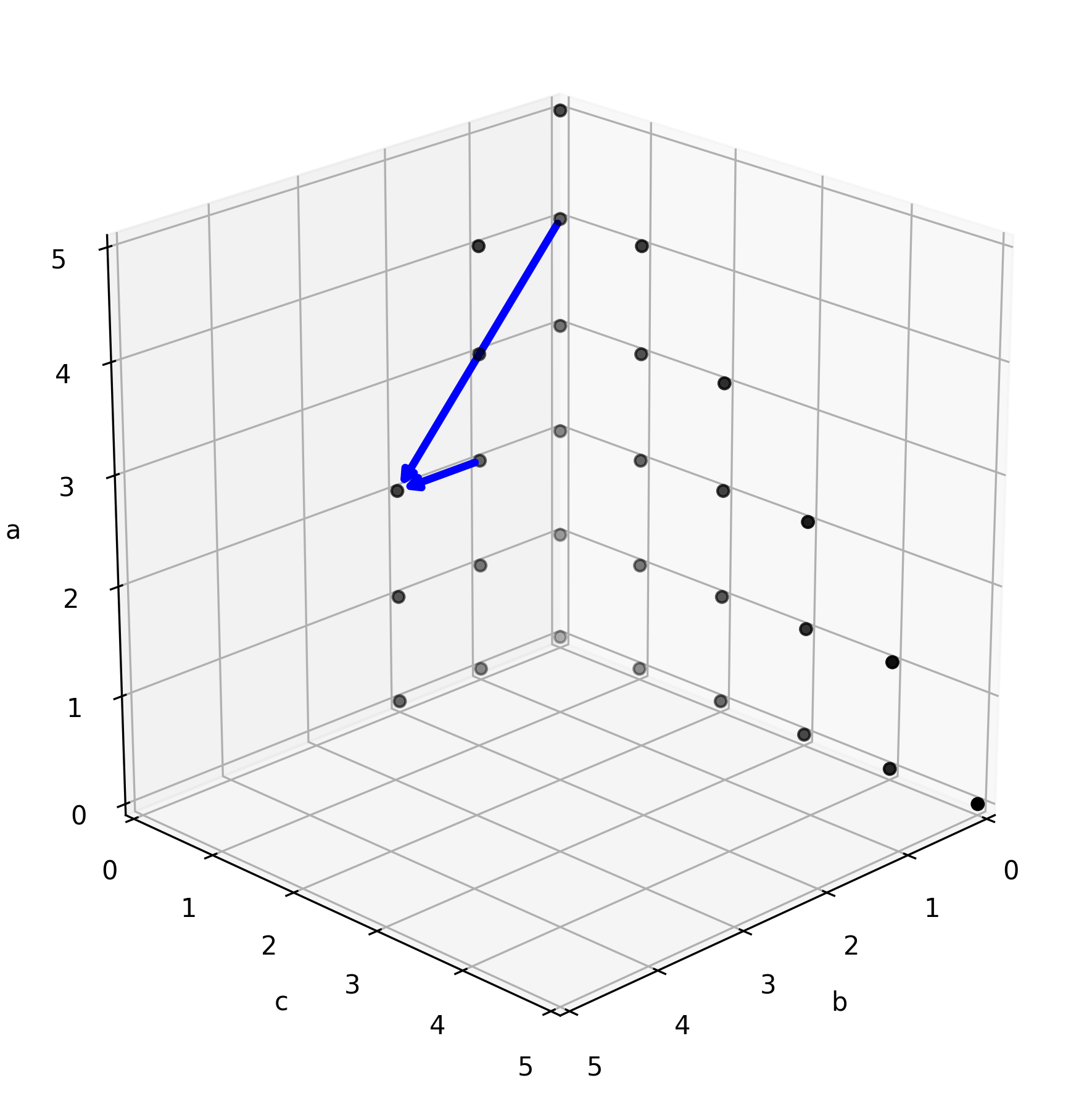}
\includegraphics[width=4cm]{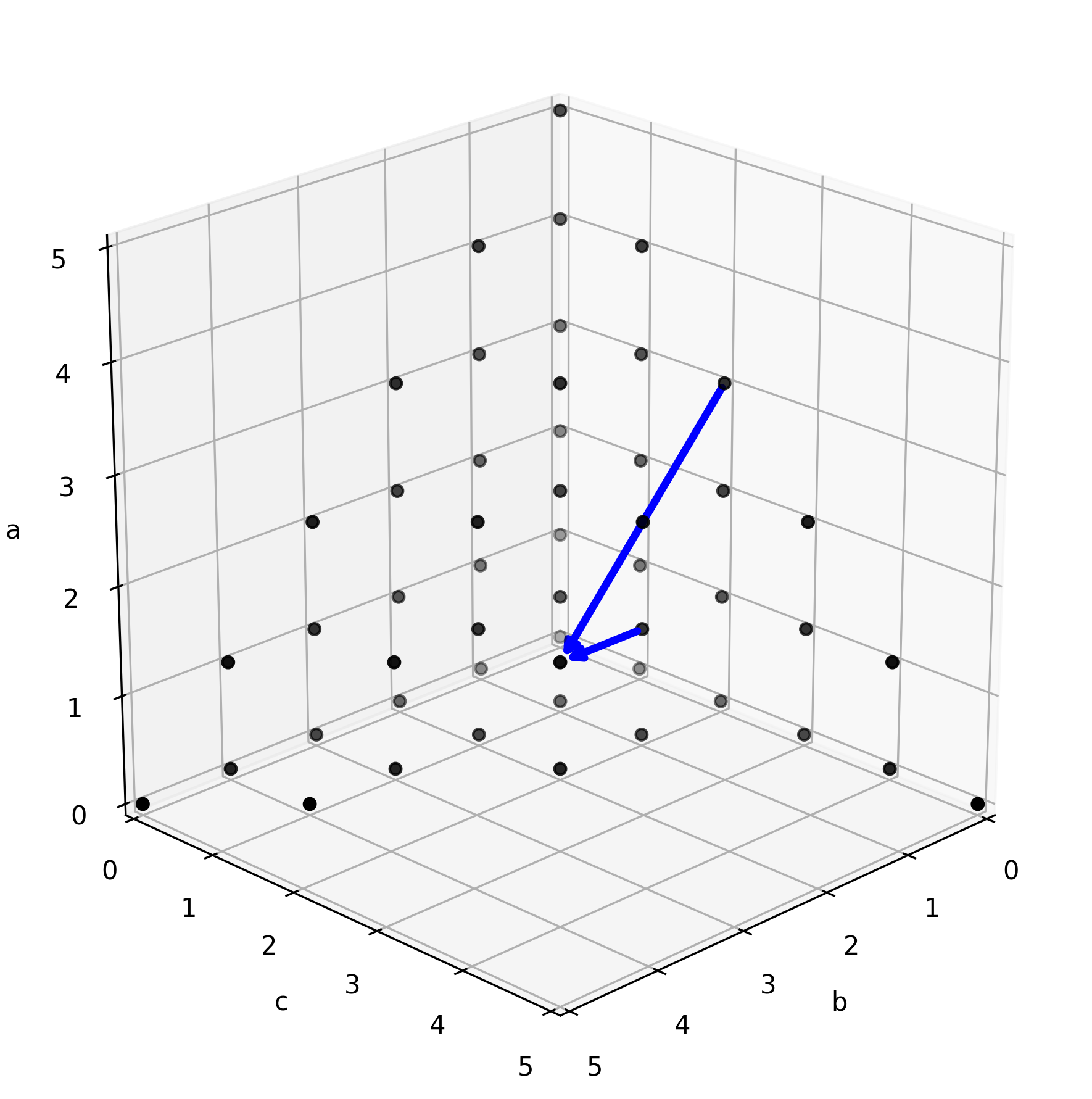}
\includegraphics[width=4cm]{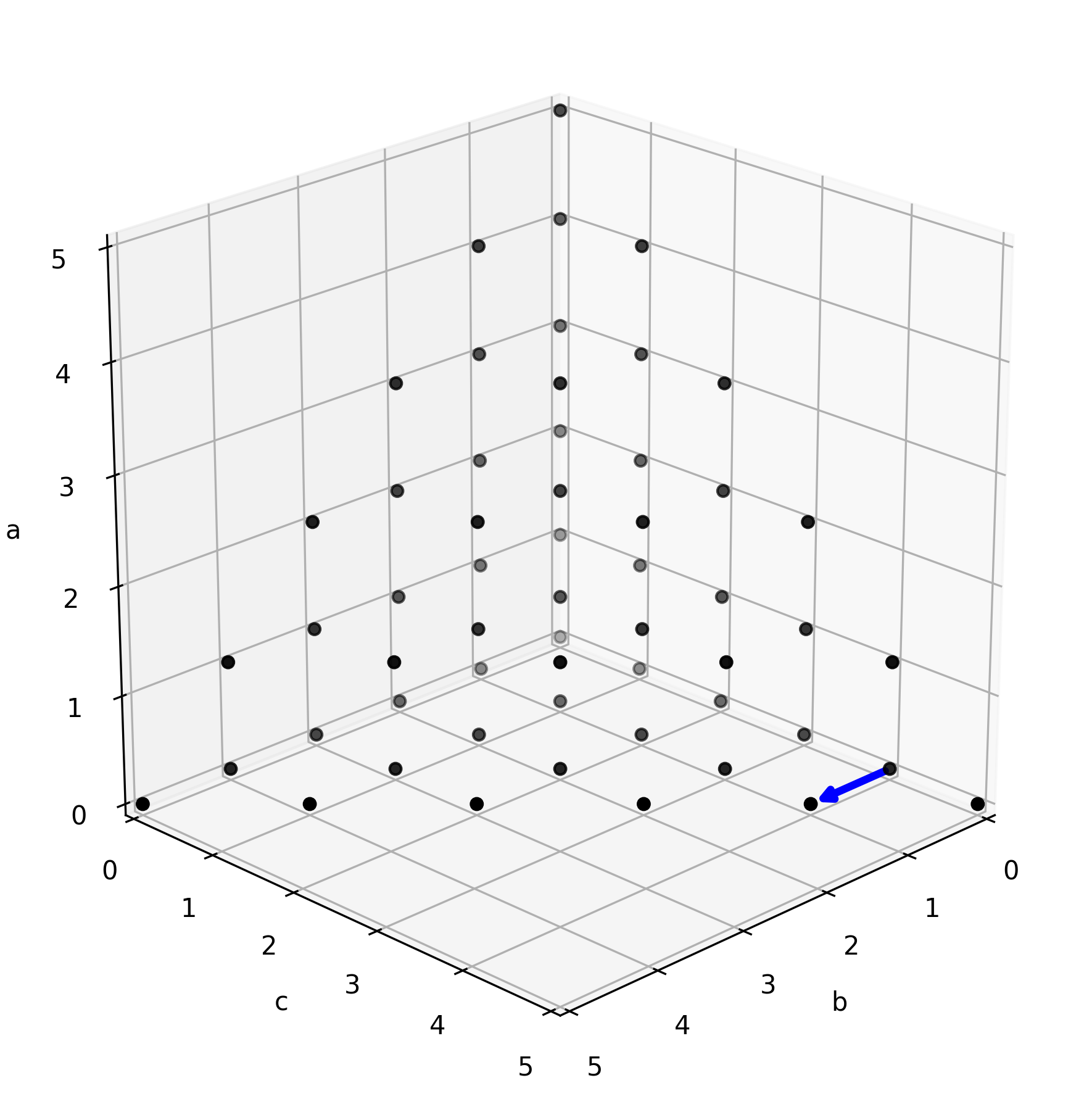}
\caption{Computing $\xi_{b,c,a}$ using \cref{algo:integral_highOrderShapeFncs_laplace_xiTable}. The blue arrows indicate terms in the recurrence relations used in each step, here for the Laplace single layer potential with $p_s=5$.
Step 1 (top left): compute seed term $\xi_{0,0,0}$, step 2 (top center): compute $\xi_{0,0,a}$ using recurrences along the $a$-axis, step 3 (top right): compute $\xi_{0,c,a}$ on the $b=0$ plane using recurrences in the $c$-directions, step 4-6 (bottom left to right): compute $\xi_{b,c,a}$ for $b>0$ using recurrences in the $b$-direction.}
\label{fig:recursionAnimated}
\end{figure}

Note that we also need to compute the singularity term associated with $\m_{b,c,a}$:
\begin{equation}\begin{aligned}
T_{b,c,a} \equiv \lim_{\varepsilon\to0} \oint_{C_{\varepsilon}} \m_{b,c,a} \cdot \n_{C_{\varepsilon}} dC_{\varepsilon},
\end{aligned}\label{eq:def_Tbca}\end{equation}
if $\check{\mathbf{r}}_p$ falls onto the element and if the integrand has a singularity on the element,
with $C_{\varepsilon}$ a circle of radius $\varepsilon$ centered at the singularity and $\n_{C_{\varepsilon}}$ the inward unit normal vector along the circle. 
This singularity term also obeys the recursions:
\begin{equation}\begin{aligned}
\!\!T_{b+1,c,a}\!\! = \!\! \frac{-b}{a+1} T_{b-1,c,a+2} + x_p T_{b,c,a}, \quad 
T_{b,c+1,a} \!\!=\!\!  \frac{-c}{a+1} T_{b,c-1,a+2} + y_p T_{b,c,a}, \\
\end{aligned}\label{eq:rec_Tsingular_laplace_bcaxis}\end{equation}
and, if the fields \cref{eq:m_0a0_laplace} are used for $b=c=0$, the following are initial values:
\begin{equation}\begin{aligned}
T_{0,0,a} 
=  \lim_{\varepsilon\to0} \oint_{C_{\varepsilon}} \m_{0,0,a} \cdot \n_{C_{\varepsilon}} dC_{\varepsilon} 
=  \lim_{\varepsilon\to0} \oint_{C_{\varepsilon}}  \frac{r^{a+1} \boldsymbol{\rho} \cdot \n_{C_{\varepsilon}}  }{(a+1)4\pi\rho^2} dC_{\varepsilon} 
= - \frac{|h|^{a+1}}{2(a+1)}. \\
\end{aligned}\label{eq:rec_Tsingular_laplace_aaxis}\end{equation}
Thus we can evaluate all terms needed for the single layer potential with an arbitrary polynomial shape function of order $p_s$ by running the above recursions. 
The procedures for computing $\kappa_{b,c,a}$ and $\xi_{b,c,a}$ are given in \cref{algo:integral_highOrderShapeFncs_laplace_etaTable} and 
\cref{algo:integral_highOrderShapeFncs_laplace_xiTable}. 
The values $\xi_{b,c,0}$ are the line integrals needed to compute the layer potential. 
The steps for computing the $\xi_{b,c,a}$ coefficients using \cref{algo:integral_highOrderShapeFncs_laplace_xiTable} are illustrated in \cref{fig:recursionAnimated} in the $(b,c,a)$ index space for the Laplace single layer potential. 
By taking derivatives and second derivatives with respect to $z$ and/or $\r_p$, similar recurrences are derived for the double layer, adjoint double layer, and hypersingular potentials. The same scheme applies for the evaluation of these layer potentials, also for the Helmholtz kernel as shown later. 
\begin{algorithm}[htbp]
\caption{Compute $\kappa_{b,c,a}$ for $b+c+a\le p_s$}\label{algo:integral_highOrderShapeFncs_laplace_etaTable}
\begin{algorithmic}
\STATE{1. Compute all $\kappa_{0,0,a}$ terms by running the recursion \cref{eq:rec_of_kappa_00a_01a_10a_11a_laplace} on the $a$-axis.}
\STATE{2. Compute all $\kappa_{0,1,a}$ terms using \cref{eq:rec_of_kappa_00a_01a_10a_11a_laplace}.}
\STATE{3. Compute all $\kappa_{0,c,a}$ terms by running recursion \cref{eq:rec_of_kappa_bplus_cplus_laplace} on the $b=0$ plane for $c=2,3,...,p_s$.}
\STATE{4. Compute all $\kappa_{1,0,a}$ terms using \cref{eq:rec_of_kappa_00a_01a_10a_11a_laplace}.}
\STATE{5. Compute all $\kappa_{b,0,a}$ terms by running recursion \cref{eq:rec_of_kappa_bplus_cplus_laplace} on the $c=0$ plane for $b=2,3,...,p_s$.}
\STATE{6. Compute all remaining $\kappa_{b+1,c+1,a}$ terms for  each $c$-plane with $c=0,1,...,p_s-1$, for $b=0,1,...,p_s-c-1$, 
using recursion \cref{eq:rec_of_kappa_bplus_cplus_laplace}. }
\end{algorithmic}
\end{algorithm}
\begin{algorithm}[htbp]
\caption{Compute $\xi_{b,c,a}$ for $b+c+a\le p_s$}\label{algo:integral_highOrderShapeFncs_laplace_xiTable}
\begin{algorithmic}
\STATE{1. Compute the seed term $\xi_{0,0,0}$ using the expression \cref{eq:rec_of_xi_00a_laplace}.}
\STATE{2. Compute all $\xi_{0,0,a}$ terms along the $a$-axis using the expression \cref{eq:rec_of_xi_00a_laplace}.}
\STATE{3. Compute all $\xi_{0,c,a}$ terms in the $b=0$ plane for $c=1,2,...,p_s$, using \cref{eq:rec_xi_bca_laplace_c}.}
\STATE{4. Compute all remaining $\xi_{b,c,a}$ terms for each $c$-plane with $c=0,1,...,p_s$, for $b=1,2,...,p_s-c$, 
using recursion \cref{eq:rec_xi_bca_laplace_b}. }
\end{algorithmic}
\end{algorithm}

Here is our main theorem:
\begin{theorem}[Layer potential integral by recursion via auxiliary vector fields - Laplace kernel]\label{thm:integral_by_recursion_laplace}
  \cref{algo:integral_highOrderShapeFncs_laplace_etaTable} and \cref{algo:integral_highOrderShapeFncs_laplace_xiTable} allow evaluation of the indefinite integral $\xi_{b,c,0} = \int \m_{b,c}\cdot\n^jdx^j$ needed for the computation of the single layer potential weighted by shape functions $x^b y^c$ based on \cref{eq:layerPotentials}.
\end{theorem}
\begin{proof}
The theorem follows from the initial conditions \cref{eq:rec_of_kappa_00a_01a_10a_11a_laplace}, \cref{eq:rec_of_xi_00a_laplace}, the recurrence relations \cref{eq:rec_of_kappa_bplus_cplus_laplace}, \cref{eq:rec_xi_bca_laplace_b} and \cref{eq:rec_xi_bca_laplace_c} due to \cref{lemma:recurrence_b_c_laplace}, and induction.
\end{proof}
Alternatively, the recursions can be executed such that $\xi_{b,c,0}$ are obtained in ascending order of $b+c$. This may be useful if application using RIPE has a mechanism to evaluate and accept results with dynamically truncated orders. 
\subsubsection{Double layer potential}
By taking the derivative of \cref{eq:rec_xi_bca_laplace_b} and \cref{eq:rec_xi_bca_laplace_c} with respect to $z$, we obtain recursions for computing the double layer potential:
\begin{equation}\begin{aligned}
  \xi_{b+1,c,a}^{(z)} &=  \frac{-b}{a+1} \xi_{b-1,c,a+2}^{(z)} + x_p \xi_{b,c,a}^{(z)} + \frac{\ihat \cdot \n^j }{4\pi} z_d \int x^{b} y^c r^{a-1} dx^j,\\
\quad \xi_{b,c+1,a}^{(z)} &= \frac{-c}{a+1} \xi_{b,c-1,a+2}^{(z)} + y_p \xi_{b,c,a}^{(z)} + \frac{\jhat \cdot \n^j }{4\pi} z_d \int x^{b} y^c r^{a-1} dx^j,\quad 
\xi_{0,0,a}^{(z)} = \frac{z_p^j k_{a-1}}{4\pi} , \\
\end{aligned}\label{eq:rec_xi_bca_laplace_dlp}\end{equation}
with $\xi_{b,c,a}^{(z)} \equiv   \partial_z \xi_{b,c,a}$.
For the singularity terms $T_{b,c,a}^{(z)} \equiv \partial_z T_{b,c,a}$, from \cref{eq:rec_Tsingular_laplace_bcaxis} we have:
\begin{equation}\begin{aligned}
&T_{b+1,c,a}^{(z)} =  x_p T_{b,c,a}^{(z)} - \frac{b}{a+1} T_{b-1,c,a+2}^{(z)} , \quad
T_{b,c+1,a}^{(z)} = y_p T_{b,c,a}^{(z)} - \frac{c}{a+1} T_{b,c-1,a+2}^{(z)} ,\\
&T_{0,0,a}^{(z)} =  \frac{z_p^j |h|^{a-1}}{2} \text{\ \ if $|h|\neq0$ else $0$}.\\ 
\end{aligned}\label{eq:rec_Tsingular_laplace_dlp}\end{equation}

\subsubsection{Adjoint double layer potential recursions}
\begin{equation}\begin{aligned}
\xi_{b+1,c,a}^{(x_p^j)}
  &= - \frac{b}{a+1} \xi_{b-1,c,a+2}^{(x_p^j)} + \alpha_x \xi_{b,c,a} + x_p \xi_{b,c,a}^{(x_p^j)} + \frac{\ihat \cdot \n^j }{(a+1)4\pi}   \partial_{x_p^j} \kappa_{b,c,a+1},\\
\xi_{b+1,c,a}^{(y_p^j)} &= - \frac{b}{a+1} \xi_{b-1,c,a+2}^{(y_p^j)} + \beta_x \xi_{b,c,a} + x_p \xi_{b,c,a}^{(y_p^j)} - \frac{\ihat \cdot \n^j }{4\pi} y_d^j \kappa_{b,c,a-1},\\
\end{aligned}\label{eq:rec_xi_bca_laplace_adlp_baxis}\end{equation}
\begin{equation}\begin{aligned}
\xi_{b,c+1,a}^{(x_p^j)} 
  &= - \frac{c}{a+1} \xi_{b,c-1,a+2}^{(x_p^j)} + \alpha_y \xi_{b,c,a}  + y_p \xi_{b,c,a}^{(x_p^j)} + \frac{\jhat \cdot \n^j }{(a+1)4\pi} \partial_{x_p^j}  \kappa_{b,c,a+1},\\
\xi_{b,c+1,a}^{(y_p^j)} &= - \frac{c}{a+1} \xi_{b,c-1,a+2}^{(y_p^j)} + \beta_y \xi_{b,c,a}  + y_p \xi_{b,c,a}^{(y_p^j)} - \frac{\jhat \cdot \n^j }{4\pi} y_d^j \kappa_{b,c,a-1},\\
\end{aligned}\label{eq:rec_xi_bca_laplace_adlp_caxis}\end{equation}
\begin{equation}\begin{aligned}
\xi_{0,0,a}^{(x_p^j)} &= \frac{-y_p^j}{(a+1)4\pi}  \frac{r^{a+1}}{\rho^2} ,\quad 
\xi_{0,0,a}^{(y_p^j)} = \frac{1}{4\pi} \left( i_{a-1} - \frac{x_d^j r^{a+1}}{(a+1)\rho^2}\right) ,\quad \xi_{b,c,a}^{(z_p^j)} = -\xi_{b,c,a}^{(z)},
\end{aligned}\label{eq:rec_xi_bca_laplace_adlp_axis}\end{equation}
with $\xi_{b,c,a}^{(x_p^j)} \equiv   \partial_{x_p^j} \xi_{b,c,a}$ and  $\xi_{b,c,a}^{(y_p^j)} \equiv   \partial_{y_p^j} \xi_{b,c,a}$.
For the singularity terms $T_{b,c,a}^{(x_p)} \equiv \partial_{x_p} T_{b,c,a}$ and $T_{b,c,a}^{(y_p)} \equiv \partial_{y_p} T_{b,c,a}$ we have:
\begin{equation}\begin{aligned}
T_{b+1,c,a}^{(x_p)} &= x_p T_{b,c,a}^{(x_p)} - \frac{b}{a+1} T_{b-1,c,a+2}^{(x_p)} + T_{b,c,a}, \quad
T_{b+1,c,a}^{(y_p)}  = x_p T_{b,c,a}^{(y_p)} - \frac{b}{a+1} T_{b-1,c,a+2}^{(y_p)}, \\
\end{aligned}\label{eq:rec_Tsingular_laplace_adlp_baxis}\end{equation}
\begin{equation}\begin{aligned}
T_{b,c+1,a}^{(x_p)} &= y_p T_{b,c,a}^{(x_p)} - \frac{c}{a+1} T_{b,c-1,a+2}^{(x_p)}, \quad
T_{b,c+1,a}^{(y_p)}  = y_p T_{b,c,a}^{(y_p)} - \frac{c}{a+1} T_{b,c-1,a+2}^{(y_p)} + T_{b,c,a}, \\
\end{aligned}\label{eq:rec_Tsingular_laplace_adlp_caxis}\end{equation}
\begin{equation}\begin{aligned}
T_{0,0,a}^{(x_p)} &= T_{0,0,a}^{(y_p)}  = 0, \quad T_{b,c,a}^{(z_p^j)} = -T_{b,c,a}^{(z)}.\\  
\end{aligned}\label{eq:rec_Tsingular_laplace_adlp_axis}\end{equation}

\subsubsection{Hypersingular Laplace potential}
For this potential we obtain:
\begin{equation}\begin{aligned}
\xi_{b+1,c,a}^{(z,x_p^j)} =& - \frac{b}{a+1} \xi_{b-1,c,a+2}^{(z,x_p^j)} + \alpha_x \xi_{b,c,a}^{(z)} + x_p \xi_{b,c,a}^{(z,x_p^j)} + \frac{\ihat \cdot \n^j }{(a+1)4\pi} \kappa_{b,c,a+1}^{(z,x_p^j)},\\
\xi_{b+1,c,a}^{(z,y_p^j)} =& - \frac{b}{a+1} \xi_{b-1,c,a+2}^{(z,y_p^j)} + \beta_x \xi_{b,c,a}^{(z)} + x_p \xi_{b,c,a}^{(z,y_p^j)} - \frac{\ihat \cdot \n^j (a-1)}{4\pi} y_d^j z_d \int x^{b} y^c r^{a-3} dx^j,\\
\xi_{b+1,c,a}^{(z,z_p^j)} =&  -\frac{b}{a+1} \xi_{b-1,c,a+2}^{(z,z_p^j)} + x_p \xi_{b,c,a}^{(z,z_p^j)} - \frac{\ihat \cdot \n^j }{4\pi} \left(  \kappa_{b,c,a-1} + z_d^2 (a-1)  \kappa_{b,c,a-3}  \right),\\
\end{aligned}\label{eq:rec_xi_bca_laplace_hsp_baxis}\end{equation}
\begin{equation}\begin{aligned}
\xi_{b,c+1,a}^{(z,x_p^j)} =& - \frac{c}{a+1} \xi_{b,c-1,a+2}^{(z,x_p^j)} + \alpha_y \xi_{b,c,a}^{(z)}  + y_p \xi_{b,c,a}^{(z,x_p^j)} + \frac{\jhat \cdot \n^j }{(a+1)4\pi} \kappa_{b,c,a+1}^{(z,x_p^j)},\\
\xi_{b,c+1,a}^{(z,y_p^j)} =& - \frac{c}{a+1} \xi_{b,c-1,a+2}^{(z,y_p^j)} + \beta_y \xi_{b,c,a}^{(z)}  + y_p \xi_{b,c,a}^{(z,y_p^j)} - \frac{\jhat \cdot \n^j (a-1)}{4\pi} y_d^j z_d \int x^{b} y^c r^{a-3} dx^j,\\
\xi_{b,c+1,a}^{(z,z_p^j)} =& -\frac{c}{a+1} \xi_{b,c-1,a+2}^{(z,z_p^j)} + y_p \xi_{b,c,a}^{(z,z_p^j)} - \frac{\jhat \cdot \n^j }{4\pi}  \left(  \kappa_{b,c,a-1} + z_d^2 (a-1)  \kappa_{b,c,a-3}  \right),\\
\end{aligned}\label{eq:rec_xi_bca_laplace_hsp_caxis}\end{equation}
\begin{equation}\begin{aligned}
\xi_{0,0,a}^{(z,x_p^j)} =&   \frac{1}{4\pi} \partial_{x_p^j}(z_p^j k_{a-1}) 
  =  \frac{1}{4\pi} z_p^j \partial_{x_p^j} k_{a-1} = \frac{1}{4\pi} z_p^j y_p^j \frac{r^{a-1}}{\rho^2}\\
\xi_{0,0,a}^{(z,y_p^j)} =&   \frac{1}{4\pi} \partial_{y_p^j}(z_p^j k_{a-1}) 
  =  \frac{1}{4\pi} z_p^j \partial_{y_p^j} k_{a-1} = \frac{1}{4\pi} z_p^j (x_d^j \frac{r^{a-1}}{\rho^2} - (a-1) i_{a-3})\\
\xi_{0,0,a}^{(z,z_p^j)} =&   \frac{1}{4\pi} \partial_{z_p^j}(z_p^j k_{a-1}) 
  =  \frac{k_{a-1}}{4\pi} + \frac{a-1}{4\pi} (z_p^j)^2 k_{a-3} 
  =  \frac{1}{4\pi} (a k_{a-1} + (a-1) y_p^j i_{a-3}), \\
\end{aligned}\label{eq:rec_xi_bca_laplace_hsp_aaxis}\end{equation}
with $\xi_{b,c,a}^{(z,x_p^j)} \equiv   \partial^2_{z,x_p^j} \xi_{b,c,a}$ and  $\xi_{b,c,a}^{(z,y_p^j)} \equiv   \partial^2_{z,y_p^j} \xi_{b,c,a}$.

For the singularity terms $T_{b,c,a}^{(z,x_p)} \equiv \partial^2_{z,x_p} T_{b,c,a}$ and $T_{b,c,a}^{(z,y_p)} \equiv \partial^2_{z,y_p} T_{b,c,a}$ we have:
\begin{equation}\begin{aligned}
T_{b+1,c,a}^{(z,x_p)} &= - \frac{b}{a+1} T_{b-1,c,a+2}^{(z,x_p)} + x_p T_{b,c,a}^{(z,x_p)} + T_{b,c,a}^{(z)}, \\
T_{b+1,c,a}^{(z,y_p)} &= - \frac{b}{a+1} T_{b-1,c,a+2}^{(z,y_p)} + x_p T_{b,c,a}^{(z,y_p)}, \quad
T_{b+1,c,a}^{(z,z_p)} = - \frac{b}{a+1} T_{b-1,c,a+2}^{(z,z_p)} + x_p T_{b,c,a}^{(z,z_p)}, \\
\end{aligned}\label{eq:rec_Tsingular_laplace_hsp_baxis}\end{equation}
\begin{equation}\begin{aligned}
T_{b,c+1,a}^{(z,y_p)} &= - \frac{c}{a+1} T_{b,c-1,a+2}^{(z,y_p)} + y_p T_{b,c,a}^{(z,y_p)} + T_{b,c,a}^{(z)}, \\
T_{b,c+1,a}^{(z,x_p)} &= - \frac{c}{a+1} T_{b,c-1,a+2}^{(z,x_p)} + y_p T_{b,c,a}^{(z,x_p)}, \quad
T_{b,c+1,a}^{(z,z_p)} = - \frac{c}{a+1} T_{b,c-1,a+2}^{(z,z_p)} + y_p T_{b,c,a}^{(z,z_p)}, \\
\end{aligned}\label{eq:rec_Tsingular_laplace_hsp_caxis}\end{equation}
\begin{equation}\begin{aligned}
T_{0,0,a}^{(z,x_p^j)} &= T_{0,0,a}^{(z,y_p^j)}  = 0, \quad T_{0,0,a}^{(z,z_p^j)} =  \frac{|h|^{a-1} + (z_p^j)^2 (a-1) |h|^{a-3}}{2} \text{\ \ if $|h|\neq0$ else $0$}.\\  
\end{aligned}\label{eq:rec_Tsingular_laplace_hsp_aaxis}\end{equation}

\subsection{Helmholtz kernel with polynomials of arbitrary order}

\subsubsection{Single layer potential}
We want to find vector fields in the form:
\begin{equation}\begin{aligned}
\nabla_s \cdot \mathbf{m}_{b,c} \equiv \frac{ x^b y^c e^{i k r}}{4 \pi   r}.\\
\end{aligned}\label{eq:divs_m_bc_helmholtz}\end{equation}

For this, the following auxiliary vector fields are introduced:
\begin{equation}\begin{aligned}
\nabla_s \cdot \mathbf{m}_{b,c,a} \equiv \frac{ x^b y^c r^{a} e^{i k r}}{4 \pi   r}.\\
\end{aligned}\label{eq:divs_m_bca_helmholtz}\end{equation}

\begin{lemma}[Recurrence relation for auxiliary vector fields - Helmholtz kernel]\label{lemma:rec_b_c_helmholtz}
  If vector fields $\m_{b+1,c,a}$, $\m_{b,c,a}$, $\m_{b-1,c,a+1}$, and $\m_{b,c,a-1}$ satisfy \cref{eq:divs_m_bca_helmholtz} and $a\ge1$, then the vector field $\m_{b+1,c,a-1}$ given by the following equation:
   \begin{equation}\begin{aligned}
\frac{a}{ik} \m_{b+1,c,a-1}= 
  x_p \m_{b,c,a}  - \m_{b+1,c,a} 
- \frac{b}{ik} \m_{b-1,c,a+1}  
+ \frac{ax_p}{ik} \m_{b,c,a-1} 
+ \frac{x^{b} y^{c} r^a e^{i k r}}{4 \pi ik}\ihat\\
\end{aligned}\label{eq:rec_m_bca_helmholtz}\end{equation}
 also satisfies \cref{eq:divs_m_bca_helmholtz} for the index tuple $(b+1,c,a-1)$.
\end{lemma}
\begin{proof}
 If vector fields $\m_{b,c,a}$, $\m_{b+1,c,a}$, $\m_{b-1,c,a+1}$ and $\m_{b,c,a-1}$ satisfy \cref{eq:divs_m_bca_helmholtz}, the surface divergence of the right hand side of \cref{eq:rec_m_bca_helmholtz} yields:
$\frac{a x^{b+1} y^c r^{a-1}e^{ikr}}{4\pi ikr}$,
hence, $\nabla_s \cdot \mathbf{m}_{b+1,c,a-1}=\frac{ x^{b+1} y^c r^{a-1}e^{ikr}}{4\pi r}$, which means that $\m_{b+1,c,a-1}$ satisfies \cref{eq:divs_m_bca_helmholtz} for the index tuple $(b+1,c,a-1)$.
\end{proof}
This yields the following recurrence relation and its counterpart due to symmetry:
\begin{equation}\begin{aligned}
\xi_{b+1,c,a} = x_p \xi_{b,c,a} -\frac{a}{ik} \xi_{b+1,c,a-1} - \frac{b}{ik} \xi_{b-1,c,a+1}  + \frac{a x_p}{ik} \xi_{b,c,a-1} +\frac{\n^j \cdot \ihat}{4\pi ik} \theta_{b,c,a},\\
\xi_{b,c+1,a} = y_p \xi_{b,c,a} - \frac{a}{ik} \xi_{b,c+1,a-1} - \frac{c}{ik} \xi_{b,c-1,a+1}  + \frac{a y_p}{ik} \xi_{b,c,a-1} +\frac{\n^j \cdot \jhat}{4\pi ik} \theta_{b,c,a},\\
\end{aligned}\label{eq:rec_xi_bca_helmholtz_bc}\end{equation}
where we have defined:
\begin{equation}\begin{aligned}
\xi_{b,c,a} \equiv \int \m_{b,c,a} \cdot \n^j dx^j, \quad \theta_{b,c,a} \equiv \int x^b y^c r^a e^{ikr} dx^j.
\end{aligned}\end{equation}
The challenge is to evaluate $\theta_{b,c,a}$.
Assuming $k|r-r_0| \le 1$ with $r_0 \equiv |\rp-\v^j|$, which is usually satisfied in the BEM as the mesh size is typically set to be smaller than $\lambda/6$ with $\lambda$ the wavelength, we expand the $e^{ikr}$ term as Taylor series:
\begin{equation}\begin{aligned}
\theta_{b,c,a} &= \int x^b y^c r^a  e^{ikr} dx^j
= e^{ikr_0} \int x^b y^c r^a   e^{ik(r-r_0)} dx^j\\
&\approx e^{ikr_0} \sum_{l=0}^{p_e-1} \frac{(ik)^l}{l!} \int x^b y^c r^a (r-r_0)^l dx^j 
 = e^{ikr_0} \sum_{l=0}^{p_e-1} A_l^{(p_e)} \int x^b y^c r^{a+l} dx^j\\
&= e^{ikr_0} \sum_{l=0}^{p_e-1} A_l^{(p_e)} \kappa_{b,c,a+l}, \quad A_l^{(p)} \equiv \frac{(ik)^l}{l!}\sum_{m=0}^{p-l-1} \frac{(-ikr_0)^m}{m!}.\\
\end{aligned}\label{eq:convolution_to_obtain_theta_bca}\end{equation}
For later use, we also define:
\begin{equation}\begin{aligned}
S_s \equiv \frac{e^{ikr_0} }{4\pi}  \sum_{l=0}^{p_e-1}A_l^{(p_e)} k_{l+s}, \quad U_s \equiv \frac{e^{ikr_0} }{4\pi}  \sum_{l=0}^{p_e-1}A_l^{(p_e)} i_{l+s}.
\end{aligned}\end{equation}
The following derivatives of these functions are also used later for the computation of derivatives of the single layer potential:
\begin{equation}\begin{aligned}
\partial_{x_p^j} S_{a+1} 
  & = \frac{e^{ikr_0} }{4\pi}  \sum_{l=0}^{p_e-1}A_l^{(p_e)} \partial_{x_p^j} k_{l+a+1} 
  = \frac{e^{ikr_0} }{4\pi}  \sum_{l=0}^{p_e-1}A_l^{(p_e)} y_p^j \frac{r^{l+a+1}}{\rho^2} 
  \approx \frac{y_p^j r^{a+1} e^{ikr} }{4\pi \rho^2},\\
  \end{aligned}\end{equation}
\begin{equation}\begin{aligned}
\partial_{y_p^j} S_{a+1} 
  & = \frac{e^{ikr_0} }{4\pi}  \sum_{l=0}^{p_e-1}A_l^{(p_e)} \left(x_d^j \frac{r^{l+a+1}}{\rho^2} - (l+a+1) i_{l+a-1} \right)\\
  & \approx \frac{x_d^j r^{a+1} e^{ikr} }{4\pi \rho^2} -  \frac{e^{ikr_0} }{4\pi} \left( \sum_{l=0}^{p_e-1}A_l^{(p_e)} (l i_{l+a-1} +   (a+1)  i_{l+a-1} )\right)\\
  & \approx \frac{x_d^j r^{a+1} e^{ikr} }{4\pi \rho^2} -  \frac{e^{ikr_0} }{4\pi}  \sum_{l=0}^{p_e-1} A_l^{(p_e)} ik i_{l+a} -   (a+1) U_{a-1}\\
  & = \frac{x_d^j r^{a+1} e^{ikr} }{4\pi \rho^2} - (ik U_a + (a+1) U_{a-1}),\\
    \end{aligned}\end{equation}
\begin{equation}\begin{aligned}
\partial_z S_{a+1} 
  & = \frac{e^{ikr_0} }{4\pi}  \sum_{l=0}^{p_e-1}A_l^{(p_e)} (l+a+1) z_d k_{l+a-1}
   \approx z_d (ik  S_a +  (a+1)   S_{a-1}),\\
    \end{aligned}\end{equation}
\begin{equation}\begin{aligned}
\partial_z U_{a} 
    = \frac{e^{ikr_0} }{4\pi}  \sum_{l=0}^{p_e-1} A_l^{(p_e)} (l+a) z_d i_{l+a-2}
  \approx z_d (ik  U_{a-1} +  a   U_{a-2}).\\
\end{aligned}\end{equation}
Thus, by computing the values of $\kappa_{b,c,a}$ with maximum $a$ index of $p_s+p_e-1$, the approximation of $\theta_{b,c,a}$ can be computed.
Note that \cref{eq:convolution_to_obtain_theta_bca} is a discrete convolution along the $a$-axis in the coefficient grid, which may be accelerated via fast Fourier transform (FFT) when both $p_s$ and $p_e$ are large.

The following recurrence relation can be used to generate vector fields $\m_{0,0,a}$ which are needed to initiate the recursions in the $b$- and $c$-axis \cref{eq:rec_xi_bca_helmholtz_bc}:
\begin{equation}\begin{aligned}
\m_{0,0,a+1} = -\frac{a+1}{ik} \m_{0,0,a} + \frac{1}{ik}\boldsymbol{\chi}_{a+1}, \quad \m_{0,0,0} = \frac{1}{ik}\boldsymbol{\chi}_0, \quad  \boldsymbol{\chi}_{a} \equiv \frac{r^a e^{ikr}}{4\pi \rho^2} \boldsymbol{\rho}.
\end{aligned}\label{eq:rec_m00a_helmholtz}\end{equation}
This recurrence yields the recurrence for the $\xi_{b,c,a}$ coefficients:
\begin{equation}\begin{aligned}
\xi_{0,0,a+1} &=  -\frac{a+1}{ik} \xi_{0,0,a} + \frac{1}{ik} \int \boldsymbol{\chi}_{a+1} \cdot \n^j dx^j
= -\frac{a+1}{ik} \xi_{0,0,a} + \frac{y_p^j}{ik} \int \frac{r^{a+1}e^{ikr}}{4\pi \rho^2} dx^j\\
& \approx -\frac{a+1}{ik} \xi_{0,0,a} -\frac{e^{ikr_0} }{4\pi ik}  \sum_{l=0}^{p_e-1}A_l^{(p_e)} k_{l+a+1} 
= -\frac{a+1}{ik} \xi_{0,0,a} - \frac{1}{ik}  S_{a+1}, \\
\xi_{0,0,0} &= \frac{y_p^j}{ik} \int \frac{e^{ikr}}{4\pi \rho^2} dx^j 
\approx  \frac{e^{ikr_0} }{4\pi ik}  \sum_{l=0}^{p_e-1}A_l^{(p_e)} y_p^j \int  \frac{r^l}{ \rho^2}    dx^j
= -\frac{1}{ik} S_0.
\end{aligned}\label{eq:rec_xi_helmholtz_aaxis}\end{equation}

The singularity term also obeys the recursions:
\begin{equation}\begin{aligned}
T_{b+1,c,a} = x_p  T_{b,c,a} - \frac{a}{ik} T_{b+1,c,a-1} - \frac{b}{ik} T_{b-1,c,a+1} + \frac{a x_p}{ik} T_{b,c,a-1}, \\
T_{b,c+1,a} = y_p  T_{b,c,a} - \frac{a}{ik} T_{b,c+1,a-1} - \frac{c}{ik} T_{b,c-1,a+1} + \frac{a y_p}{ik} T_{b,c,a-1}. \\
\end{aligned}\label{eq:rec_Tsingular_helmholtz_bcaxis}\end{equation}

If the vector fields $\m_{0,0,a}$ are chosen using the recurrence \cref{eq:rec_m00a_helmholtz}, we have the following recurrence along the $a$-axis:
\begin{equation}\begin{aligned}
\! T_{0,0,a+1} &= \!\! -\frac{a+1}{ik} T_{0,0,a} + \frac{1}{ik} \lim_{\varepsilon\to0} \oint_{C_{\varepsilon}}  \boldsymbol{\chi}_{a+1} \cdot \n_{C_{\varepsilon}} dC_{\varepsilon} 
= -\frac{a+1}{ik} T_{0,0,a} -  \frac{|h|^{a+1} e^{ik|h|}}{2ik}.
\end{aligned}\label{eq:rec_Tsingular_helmholtz_aaxis}\end{equation}
Note that this equation holds for $a=-1$ and gives the seed term $T_{0,0,0} = ie^{ik|h|}/(2k)$.

\subsubsection{Double layer potential}

The derivatives of the $\xi_{0,0,a}$ are given by:
\begin{equation}\begin{aligned}
\xi_{0,0,0}^{(z)} \approx 
-z_d^j S_{-1},\quad
\xi_{0,0,a+1}^{(z)} 
 \approx -\frac{a+1}{ik} \xi_{0,0,a}^{(z)} - z_d S_a -  \frac{(a+1) z_d}{ik}  S_{a-1}. \\
\end{aligned}\label{eq:rec_xi_bca_helmholtz_dlp_aaxis}\end{equation}

Recurrences along $b$- and $c$-axis:
\begin{equation}\begin{aligned}
\xi_{b+1,c,a}^{(z)} = x_p \xi_{b,c,a}^{(z)} - \frac{a}{ik} \xi_{b+1,c,a-1}^{(z)} - \frac{b}{ik} \xi_{b-1,c,a+1}^{(z)}  + \frac{ax_p}{ik} \xi_{b,c,a-1}^{(z)} +\frac{\n^j \cdot \ihat}{4\pi ik} \theta_{b,c,a}^{(z)},\\
\xi_{b,c+1,a}^{(z)} = y_p \xi_{b,c,a}^{(z)} - \frac{a}{ik} \xi_{b,c+1,a-1}^{(z)} - \frac{c}{ik} \xi_{b,c-1,a+1}^{(z)}  + \frac{a y_p}{ik} \xi_{b,c,a-1}^{(z)} +\frac{\n^j \cdot \jhat}{4\pi ik} \theta_{b,c,a}^{(z)}.\\
\end{aligned}\label{eq:rec_xi_bca_helmholtz_dlp_bcaxis}\end{equation}
Here, the derivative $\theta_{b,c,a}^{(z)} \equiv \partial_z \theta_{b,c,a}$ can be computed as:
\begin{equation}\begin{aligned}
\theta_{b,c,a}^{(z)} &\approx  e^{ikr_0} \sum_{l=0}^{p_e-1} A_l^{(p_e)} (a+l) z_d \int x^b y^c r^{a+l-2} dx^j\\
&= e^{ikr_0} a z_d \sum_{l=0}^{p_e-1} A_l^{(p_e)} \int x^b y^c r^{a+l-2} dx^j + e^{ikr_0} \sum_{l=1}^{p_e-1} A_l^{(p_e)} l z_d \int x^b y^c r^{a+l-2} dx^j\\
&= e^{ikr_0} a z_d \sum_{l=0}^{p_e-1} A_l^{(p_e)} \int x^b y^c r^{a+l-2} dx^j + e^{ikr_0} ik z_d \sum_{l=1}^{p_e-1} A_{l-1}^{(p_e-1)} \int x^b y^c r^{a+l-2} dx^j\\
&\approx e^{ikr_0} a z_d \sum_{l=0}^{p_e-1} A_l^{(p_e)} \int x^b y^c r^{a+l-2} dx^j + e^{ikr_0} ik z_d \sum_{l=0}^{p_e-1} A_{l}^{(p_e)} \int x^b y^c r^{a+l-1} dx^j\\
& = z_d(a  \theta_{b,c,a-2} + ik \theta_{b,c,a-1}).
\end{aligned}\end{equation}

If the vector fields $\m_{0,0,a}$ are chosen using the recurrence \cref{eq:rec_m00a_helmholtz}, we have the following recurrence along the $a$-axis for the singularity term:
\begin{equation}\begin{aligned}
T_{0,0,0}^{(z)} &= -  \partial_z \frac{ e^{ik|h|}}{2ik} =  \frac{-z_d  e^{i k \left|h\right|}}{2  \left|h\right|} \quad \text{if $z\neq z_p$ else 0}, \\
T_{0,0,a+1}^{(z)} 
  &= -\frac{a+1}{ik} T_{0,0,a}^{(z)} + \frac{z_d \left(i (a+1) - k \left|h\right| \right) e^{i k \left|h\right|} \left|h\right|^{a - 1}}{2 k}.
\end{aligned}\label{eq:rec_Tsingular_helmholtz_dlp_aaxis}\end{equation}
The recursion for the singularity term along the $b$- and $c$- axis is given by:
\begin{equation}\begin{aligned}
T_{b+1,c,a}^{(z)} = x_p  T_{b,c,a}^{(z)} - \frac{a}{ik} T_{b+1,c,a-1}^{(z)} - \frac{b}{ik} T_{b-1,c,a+1}^{(z)} + \frac{a x_p}{ik} T_{b,c,a-1}^{(z)}, \\
T_{b,c+1,a}^{(z)} = y_p  T_{b,c,a}^{(z)} - \frac{a}{ik} T_{b,c+1,a-1}^{(z)} - \frac{c}{ik} T_{b,c-1,a+1}^{(z)} + \frac{a y_p}{ik} T_{b,c,a-1}^{(z)}. \\
\end{aligned}\label{eq:rec_Tsingular_helmholtz_dlp_bcaxis}\end{equation}

\subsubsection{Adjoint double layer potential}

Seed term:
\begin{equation}\begin{aligned}
\xi_{0,0,0}^{(x_p^j)} &= -\frac{y_p^j e^{ikr}}{4\pi ik \rho^2} ,\quad 
\xi_{0,0,0}^{(y_p^j)} =   U_{-1} - \frac{x_d^j e^{ikr}}{4\pi ik \rho^2}.
\end{aligned}\label{eq:rec_xi_bca_helmholtz_adlp_seed}\end{equation}
Recurrence along $a$-axis is given by:
\begin{equation}\begin{aligned}
\xi_{0,0,a+1}^{(x_p^j)} 
  & = -\frac{a+1}{ik} \xi_{0,0,a}^{(x_p^j)} - \frac{1}{ik}  S_{a+1}^{(x_p^j)}  = -\frac{a+1}{ik} \xi_{0,0,a}^{(x_p^j)} - \frac{y_p^j r^{a+1} e^{ikr}}{4\pi ik \rho^2}, \\
\xi_{0,0,a+1}^{(y_p^j)} 
  & =  -\frac{a+1}{ik} \xi_{0,0,a}^{(y_p^j)} +  U_a + \frac{a+1}{ik} U_{a-1} - \frac{x_d^j r^{a+1} e^{ikr}}{4\pi ik \rho^2} .
\end{aligned}\label{eq:rec_xi_bca_helmholtz_adlp_aaxis}\end{equation}
Recurrence along $b$- and $c$-axis are given by:
\begin{equation}\begin{aligned}
\xi_{b+1,c,a}^{(x_p^j)} = & x_p \xi_{b,c,a}^{(x_p^j)} + \alpha_x \xi_{b,c,a} - \frac{a}{ik} \xi_{b+1,c,a-1}^{(x_p^j)}
 - \frac{b}{ik} \xi_{b-1,c,a+1}^{(x_p^j)} \\ &+ \frac{ax_p}{ik} \xi_{b,c,a-1}^{(x_p^j)} + \frac{a\alpha_x}{ik} \xi_{b,c,a-1} +\frac{\n^j \cdot \ihat}{4\pi ik} \theta_{b,c,a}^{(x_p^j)},\\
\xi_{b+1,c,a}^{(y_p^j)} = & x_p \xi_{b,c,a}^{(y_p^j)} + \beta_x \xi_{b,c,a} - \frac{a}{ik} \xi_{b+1,c,a-1}^{(y_p^j)} 
 - \frac{b}{ik} \xi_{b-1,c,a+1}^{(y_p^j)} \\ &+ \frac{ax_p}{ik} \xi_{b,c,a-1}^{(y_p^j)} + \frac{a\beta_x}{ik} \xi_{b,c,a-1} +\frac{\n^j \cdot \ihat}{4\pi ik} \theta_{b,c,a}^{(y_p^j)},\\ 
\end{aligned}\label{eq:rec_xi_bca_helmholtz_adlp_baxis}\end{equation}
\begin{equation}\begin{aligned}
\xi_{b,c+1,a}^{(x_p^j)} = & y_p \xi_{b,c,a}^{(x_p^j)} + \alpha_y \xi_{b,c,a} - \frac{a}{ik} \xi_{b,c+1,a-1}^{(x_p^j)}
- \frac{c}{ik} \xi_{b,c-1,a+1}^{(x_p^j)}  \\&+ \frac{a y_p}{ik} \xi_{b,c,a-1}^{(x_p^j)}+ \frac{a \alpha_y}{ik} \xi_{b,c,a-1} +\frac{\n^j \cdot \jhat}{4\pi ik} \theta_{b,c,a}^{(x_p^j)},\\
\xi_{b,c+1,a}^{(y_p^j)} = & y_p \xi_{b,c,a}^{(y_p^j)} + \beta_y \xi_{b,c,a} - \frac{a}{ik} \xi_{b,c+1,a-1}^{(y_p^j)}
- \frac{c}{ik} \xi_{b,c-1,a+1}^{(y_p^j)}  \\&+ \frac{a y_p}{ik} \xi_{b,c,a-1}^{(y_p^j)}+ \frac{a \beta_y}{ik} \xi_{b,c,a-1} +\frac{\n^j \cdot \jhat}{4\pi ik} \theta_{b,c,a}^{(y_p^j)},\\
\end{aligned}\label{eq:rec_xi_bca_helmholtz_adlp_caxis}\end{equation}
where we used:
\begin{equation}\begin{aligned}
\theta_{b,c,a}^{(x_p^j)} &\approx -x^b y^c r^a e^{ikr} + \alpha_x b \theta_{b-1,c,a} + \alpha_y c  \theta_{b,c-1,a},\\
\theta_{b,c,a}^{(y_p^j)} &\approx -y_d^j(a \theta_{b,c,a-2} + ik \theta_{b,c,a-1}).\quad
\end{aligned}\end{equation}

Recursions for the singularity term are given by: 
\begin{equation}\begin{aligned}
T_{b+1,c,a}^{(x_p)} &= x_p  T_{b,c,a}^{(x_p)} + T_{b,c,a} - \frac{a}{ik} T_{b+1,c,a-1}^{(x_p)} - \frac{b}{ik} T_{b-1,c,a+1}^{(x_p)} + \frac{a x_p}{ik} T_{b,c,a-1}^{(x_p)} + \frac{a}{ik} T_{b,c,a-1}, \\
T_{b+1,c,a}^{(y_p)} &= x_p  T_{b,c,a}^{(y_p)}  - \frac{a}{ik} T_{b+1,c,a-1}^{(y_p)}  - \frac{b}{ik} T_{b-1,c,a+1}^{(y_p)}  + \frac{a x_p}{ik} T_{b,c,a-1}^{(y_p)} , \\
\end{aligned}\label{eq:rec_Tsingular_helmholtz_adlp_baxis}\end{equation}
\begin{equation}\begin{aligned}
T_{b,c+1,a}^{(x_p)} &= y_p  T_{b,c,a}^{(x_p)} - \frac{a}{ik} T_{b,c+1,a-1}^{(x_p)} - \frac{c}{ik} T_{b,c-1,a+1}^{(x_p)} + \frac{a y_p}{ik} T_{b,c,a-1}^{(x_p)}, \\
T_{b,c+1,a}^{(y_p)} &= y_p  T_{b,c,a}^{(y_p)} + T_{b,c,a} - \frac{a}{ik} T_{b,c+1,a-1}^{(y_p)} - \frac{c}{ik} T_{b,c-1,a+1}^{(y_p)} + \frac{a y_p}{ik} T_{b,c,a-1}^{(y_p)} + \frac{a}{ik} T_{b,c,a-1}. \\
\end{aligned}\label{eq:rec_Tsingular_helmholtz_adlp_caxis}\end{equation}
Note that $T_{0,0,a}^{(x_p)} = T_{0,0,a}^{(y_p)} = 0$,
$\xi_{b,c,a}^{(z_p^j)} = - \xi_{b,c,a}^{(z)}$, and $T_{b,c,a}^{(z_p^j)} = - T_{b,c,a}^{(z)}$.

\subsubsection{Hypersingular potential}
Seed terms:
\begin{equation}\begin{aligned}
\xi_{0,0,0}^{(z,x_p^j)} &= y_p^j  z_p^j  \frac{ e^{ikr}}{4\pi \rho^2 r} ,\quad
\xi_{0,0,0}^{(z,y_p^j)} = z_d (ik  U_{-2} -   U_{-3}) + z_p^j  \frac{ x_d^j e^{ikr}}{4\pi \rho^2 r }  ,\\
\xi_{0,0,0}^{(z,z_p^j)} &= (z_d^j S_{-1})^{(z)} = S_{-1} + z_d^j S_{-1}^{(z)}  = S_{-1} + (z_p^j)^2 (ik  S_{-2} - S_{-3}).\\
\end{aligned}\label{eq:rec_xi_bca_helmholtz_hsp_seed}\end{equation}
Recurrence along $a$-axis:
\begin{equation}\begin{aligned}
\xi_{0,0,a+1}^{(z,x_p^j)} =& -\frac{a+1}{ik} \xi_{0,0,a}^{(z,x_p^j)} + y_p^j z_p^j  \frac{\left(k r - i \left(a + 1\right) \right) r^{a-1}e^{i k r}}{4 \pi k \rho^2}, \\
\xi_{0,0,a+1}^{(z,y_p^j)}  =& -\frac{a+1}{ik} \xi_{0,0,a}^{(z,y_p^j)} +  z_d (ik  U_{a-1} +  (2a+1)   U_{a-2}  + \frac{a^2-1}{ik} U_{a-3}) \\
                            & + z_p^j \frac{x_d^j \left(k r - i \left(a + 1\right) \right) r^{a-1} e^{i k r}}{4 \pi k \rho^2}, \\
\xi_{0,0,a+1}^{(z,z_p^j)}
  =& -\frac{a+1}{ik} \xi_{0,0,a}^{(z,z_p^j)} + S_a +  (\frac{a+1}{ik}+z_d^2 ik)  S_{a-1} \\
   & + z_d^2  ((2a+1) S_{a-2}  +   \frac{a^2-1}{ik} S_{a-3}). \\
\end{aligned}\label{eq:rec_xi_bca_helmholtz_hsp_aaxis}\end{equation}
Recurrence along $b$- and $c$-axis:
\begin{equation}\begin{aligned}
\xi_{b+1,c,a}^{(z,x_p^j)} = & x_p \xi_{b,c,a}^{(z,x_p^j)} + \alpha_x \xi_{b,c,a}^{(z)} - \frac{a}{ik} \xi_{b+1,c,a-1}^{(z,x_p^j)}
 - \frac{b}{ik} \xi_{b-1,c,a+1}^{(z,x_p^j)} \\ &+ \frac{ax_p}{ik} \xi_{b,c,a-1}^{(z,x_p^j)} + \frac{a\alpha_x}{ik} \xi_{b,c,a-1}^{(z)} +\frac{\n^j \cdot \ihat}{4\pi ik} \theta_{b,c,a}^{(z,x_p^j)},\\
\xi_{b+1,c,a}^{(z,y_p^j)} = & x_p \xi_{b,c,a}^{(z,y_p^j)} + \beta_x \xi_{b,c,a}^{(z)} - \frac{a}{ik} \xi_{b+1,c,a-1}^{(z,y_p^j)} 
 - \frac{b}{ik} \xi_{b-1,c,a+1}^{(z,y_p^j)} \\ &+ \frac{ax_p}{ik} \xi_{b,c,a-1}^{(z,y_p^j)} + \frac{a\beta_x}{ik} \xi_{b,c,a-1}^{(z)} +\frac{\n^j \cdot \ihat}{4\pi ik} \theta_{b,c,a}^{(z,y_p^j)},\\ 
\xi_{b+1,c,a}^{(z,z_p^j)} =& x_p \xi_{b,c,a}^{(z,z_p^j)} - \frac{a}{ik} \xi_{b+1,c,a-1}^{(z,z_p^j)} - \frac{b}{ik} \xi_{b-1,c,a+1}^{(z,z_p^j)}  + \frac{ax_p}{ik} \xi_{b,c,a-1}^{(z,z_p^j)} +\frac{\n^j \cdot \ihat}{4\pi ik} \theta_{b,c,a}^{(z,z_p^j)},\\
\end{aligned}\label{eq:rec_xi_bca_helmholtz_hsp_baxis}\end{equation}
\begin{equation}\begin{aligned}
\xi_{b,c+1,a}^{(z,x_p^j)} = & y_p \xi_{b,c,a}^{(z,x_p^j)} + \alpha_y \xi_{b,c,a}^{(z)} - \frac{a}{ik} \xi_{b,c+1,a-1}^{(z,x_p^j)}
- \frac{c}{ik} \xi_{b,c-1,a+1}^{(z,x_p^j)}  \\&+ \frac{a y_p}{ik} \xi_{b,c,a-1}^{(z,x_p^j)}+ \frac{a \alpha_y}{ik} \xi_{b,c,a-1}^{(z)} +\frac{\n^j \cdot \jhat}{4\pi ik} \theta_{b,c,a}^{(z,x_p^j)},\\
\xi_{b,c+1,a}^{(z,y_p^j)} = & y_p \xi_{b,c,a}^{(z,y_p^j)} + \beta_y \xi_{b,c,a}^{(z)} - \frac{a}{ik} \xi_{b,c+1,a-1}^{(z,y_p^j)}
- \frac{c}{ik} \xi_{b,c-1,a+1}^{(z,y_p^j)}  \\&+ \frac{a y_p}{ik} \xi_{b,c,a-1}^{(z,y_p^j)}+ \frac{a \beta_y}{ik} \xi_{b,c,a-1}^{(z)} +\frac{\n^j \cdot \jhat}{4\pi ik} \theta_{b,c,a}^{(z,y_p^j)},\\
\xi_{b,c+1,a}^{(z,z_p^j)} =& y_p \xi_{b,c,a}^{(z,z_p^j)} - \frac{a}{ik} \xi_{b,c+1,a-1}^{(z,z_p^j)} - \frac{c}{ik} \xi_{b,c-1,a+1}^{(z,z_p^j)}  + \frac{a y_p}{ik} \xi_{b,c,a-1}^{(z,z_p^j)} +\frac{\n^j \cdot \jhat}{4\pi ik} \theta_{b,c,a}^{(z,z_p^j)}.\\
\end{aligned}\label{eq:rec_xi_bca_helmholtz_hsp_caxis}\end{equation}
Here we used:
\begin{equation}\begin{aligned}
\theta_{b,c,a}^{(z,x_p^j)} &\approx -x^b y^c (a+ikr)r^{a-2}z_d e^{ikr} + \alpha_x b \theta_{b-1,c,a}^{(z)} + \alpha_y c  \theta_{b,c-1,a}^{(z)},\\
\theta_{b,c,a}^{(z,y_p^j)} &\approx -y_d^j(a \theta_{b,c,a-2}^{(z)} + ik \theta_{b,c,a-1}^{(z)}),\quad 
\theta_{b,c,a}^{(z)} \approx z_d(a  \theta_{b,c,a-2} + ik \theta_{b,c,a-1}),\\
\theta_{b,c,a}^{(z,z_p^j)} &\approx -z_d^j(a \theta_{b,c,a-2}^{(z)} + ik \theta_{b,c,a-1}^{(z)}) - (a \theta_{b,c,a-2} + ik \theta_{b,c,a-1}).\\
\end{aligned}\end{equation}

Recursions for the singularity term:
\begin{equation}\begin{aligned}
T_{0,0,a}^{(z,x_p)} &= T_{0,0,a}^{(z,y_p)} = 0, \quad T_{0,0,0}^{(z,z_p)} = \frac{ike^{ik|h|}}{2}, \\
T_{0,0,a+1}^{(z,z_p)} &= -\frac{a+1}{ik} T_{0,0,a}^{(z,z_p)} - \frac{e^{ik|h|} |h|^{a-1} }{2k}  (-2(a+1)k|h| +i(a(a+1)-k^2 |h|^2)) ,\\
\end{aligned}\label{eq:rec_Tsingular_helmholtz_hsp_aAxis}\end{equation}
\begin{equation}\begin{aligned}
T_{b+1,c,a}^{(z,x_p)} &= x_p  T_{b,c,a}^{(z,x_p)} + T_{b,c,a}^{(z)} - \frac{b}{ik} T_{b-1,c,a+1}^{(z,x_p)} + \frac{a }{ik} (x_p T_{b,c,a-1}^{(z,x_p)} +  T_{b,c,a-1}^{(z)} - T_{b+1,c,a-1}^{(z,x_p)}), \\
T_{b+1,c,a}^{(z,y_p)} &= x_p  T_{b,c,a}^{(z,y_p)}  - \frac{a}{ik} T_{b+1,c,a-1}^{(z,y_p)}  - \frac{b}{ik} T_{b-1,c,a+1}^{(z,y_p)}  + \frac{a x_p}{ik} T_{b,c,a-1}^{(z,y_p)} , \\
T_{b+1,c,a}^{(z,z_p)} &= x_p  T_{b,c,a}^{(z,z_p)} - \frac{a}{ik} T_{b+1,c,a-1}^{(z,z_p)} - \frac{b}{ik} T_{b-1,c,a+1}^{(z,z_p)} + \frac{a x_p}{ik} T_{b,c,a-1}^{(z,z_p)}, \\
\end{aligned}\label{eq:rec_Tsingular_helmholtz_hsp_bAxis}\end{equation}
\begin{equation}\begin{aligned}
T_{b,c+1,a}^{(z,x_p)} &= y_p  T_{b,c,a}^{(z,x_p)} - \frac{a}{ik} T_{b,c+1,a-1}^{(z,x_p)} - \frac{c}{ik} T_{b,c-1,a+1}^{(z,x_p)} + \frac{a y_p}{ik} T_{b,c,a-1}^{(z,x_p)}, \\
T_{b,c+1,a}^{(z,y_p)} &= y_p  T_{b,c,a}^{(z,y_p)} + T_{b,c,a}^{(z)}  - \frac{c}{ik} T_{b,c-1,a+1}^{(z,y_p)} + \frac{a }{ik} (y_p T_{b,c,a-1}^{(z,y_p)} + T_{b,c,a-1}^{(z)} - T_{b,c+1,a-1}^{(z,y_p)}), \\
T_{b,c+1,a}^{(z,z_p)} &= y_p  T_{b,c,a}^{(z,z_p)} - \frac{a}{ik} T_{b,c+1,a-1}^{(z,z_p)} - \frac{c}{ik} T_{b,c-1,a+1}^{(z,z_p)} + \frac{a y_p}{ik} T_{b,c,a-1}^{(z,z_p)}. \\
\end{aligned}\label{eq:rec_Tsingular_helmholtz_hsp_cAxis}\end{equation}

\subsection{Complexity}
For the Laplace kernel, evaluating the recursions in \cref{algo:integral_highOrderShapeFncs_laplace_etaTable} and \cref{algo:integral_highOrderShapeFncs_laplace_xiTable} has a complexity of $O(p_s^3)$. Note that this is the complexity to obtain all terms with the factor $x^by^c$ for $b+c\le p_s$. 
The layer potential for one shape function can then be computed as a linear combination of $O(p_s^2)$ terms. 
If the layer potential was to be evaluated by Gauss-Legendre quadrature of order $b+c+p_e'$ with $p_e'$ the degrees to accommodate the non-polynomial part, i.e., the Green function part of the integrand, the complexity to evaluate one layer potential would be $O((p_s + p_e')^2 )$, assuming the use of $O((p_s+p_e')^2)$ quadrature nodes. 
For the Helmholtz case, if we use $p_e$ degrees for the Taylor series expansion of the oscillatory term $e^{ikr}$, the complexity of the recursions in the RIPE method is $O(p_s^3 p_e )$ with a simple implementation and $O(p_s^2 (p_s+p_e) \log(p_s+p_e))$ if FFT is used for the convolution. For Gauss-Legendre, the complexity is $O((p_s + p_e')^2  )$, where $p_e'$ degrees are used to accommodate the non-polynomial part. Note that $p_e$ and $p_e'$, although both parameters specify the precision of the integral evaluation, cannot be directly associated; Gauss-Legendre would require larger $p_e'$ for evaluation points in the near-field since near the singularity the integrand cannot be well approximated by low-order polynomials. On the other hand, $p_e$ in the RIPE method can be smaller in the near-field as it is the expansion order of the Taylor series. This implies that the appropriate method depends on the distance of the evaluation point from the element. 
If low-order elements with $p_s = O(1)$ are used, the complexity of the RIPE method for the Laplace and Helmholtz kernel reduces to $O(1)$ and $O(p_e)$, respectively. 
The complexities for the low-order case and general case are summarized in \cref{table:complexity}. 
\begin{table}[htbp]
\begin{center}
\begin{tabular}{@{}c@{}|@{}c@{}|@{}c@{}}
    \begin{tabular}{c}
           \\
          \\\hline \hline
       RIPE    \\\hline
       RIPE (+FFT)    \\\hline
       Gauss-Legendre \\ 
    \end{tabular} &
    \begin{tabular}{c|c}
       \multicolumn{2}{c}{Laplace}\\\hline 
       $p_s \approx 1$ & General \\\hline\hline
        $O(1)$ & $O(p_s^3)$ \\\hline
       - & - \\\hline
       $O(p_e'^2)$ & 
       $O( (p_s + p_e')^2  )$ 
    \end{tabular} &
    \begin{tabular}{c|c}
       \multicolumn{2}{c}{Helmholtz}\\\hline
       $p_s \approx 1$ & General \\\hline\hline
       $O(p_e)$ & $O(p_s^3 p_e)$ \\\hline
       - & $O(p_s^2 (p_s+p_e) \log(p_s+p_e))$ \\\hline
       $O(p_e'^2)$ & $O( (p_s + p_e')^2  )$ \\
    \end{tabular} 
\end{tabular}
\caption{Complexity for evaluating layer potentials for $p_s$ of $O(1)$ and general order, respectively. Note that the complexities  for Gauss-Legendre are for one layer potential evaluation, while for RIPE they are the complexities to computes all terms with the monomial factor $x^by^c$ for $0 \le b+c\le p_s$ in a single execution of the recursive algorithm. These terms can be reused when evaluating multiple layer potentials with different linear combinations of these monomial terms. }
\label{table:complexity}
\end{center}
\end{table}

\section{Numerical evaluation}
\subsection{Accuracy}
\subsubsection{Nearly-singular case}
\begin{figure*}[htbp]
    \centering
    \includegraphics[width=6cm]{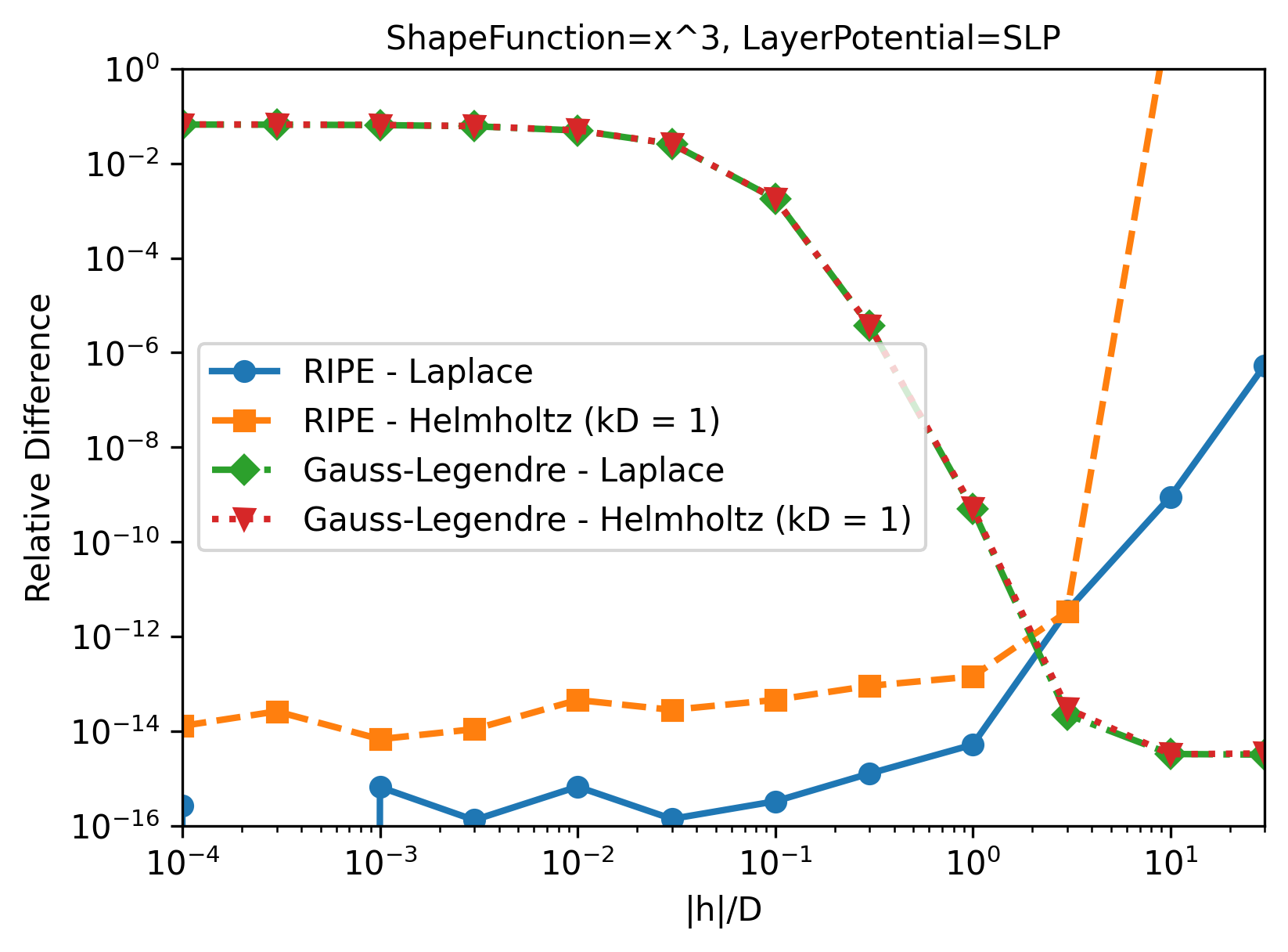}
    \includegraphics[width=6cm]{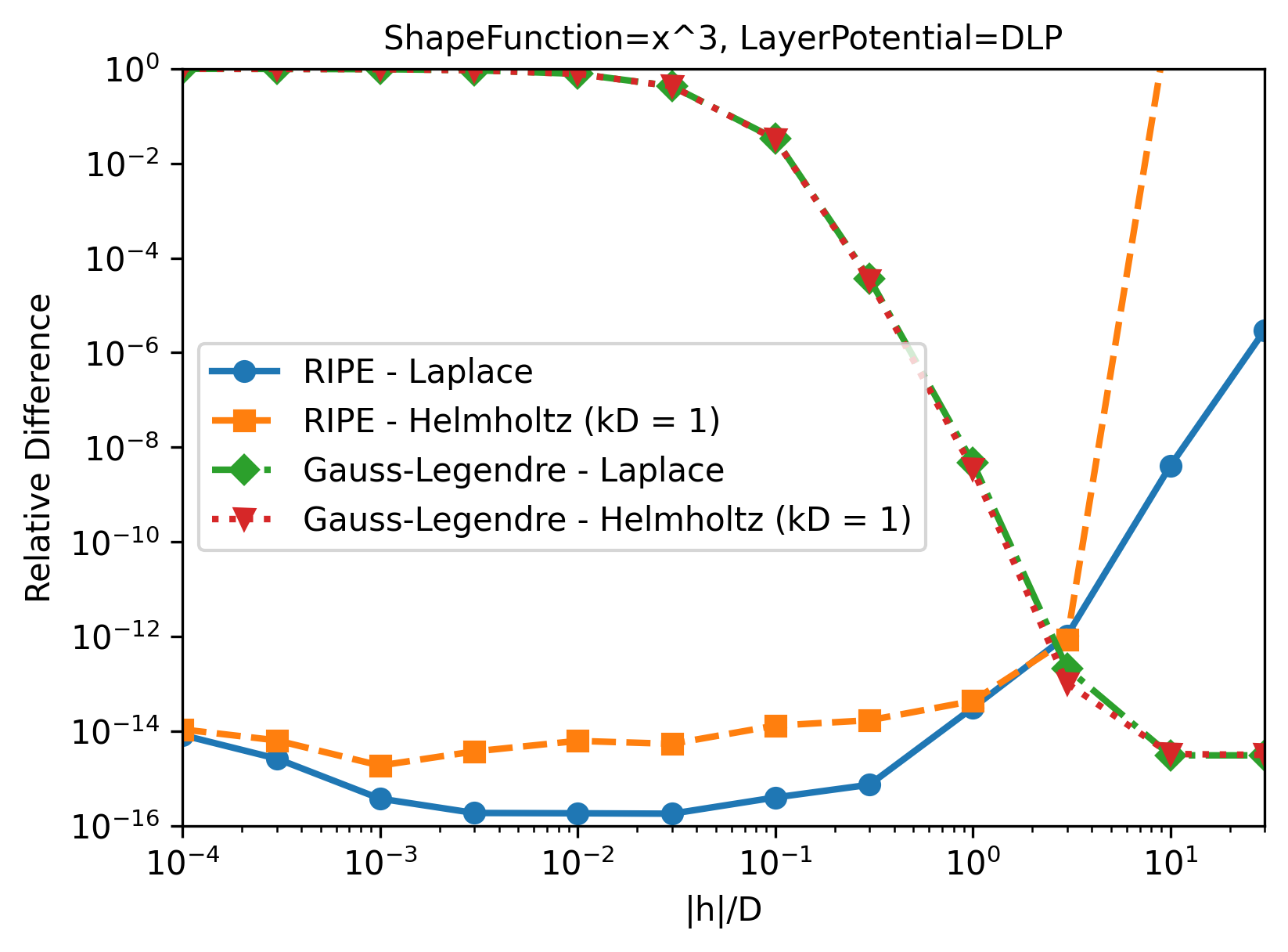}
    \includegraphics[width=6cm]{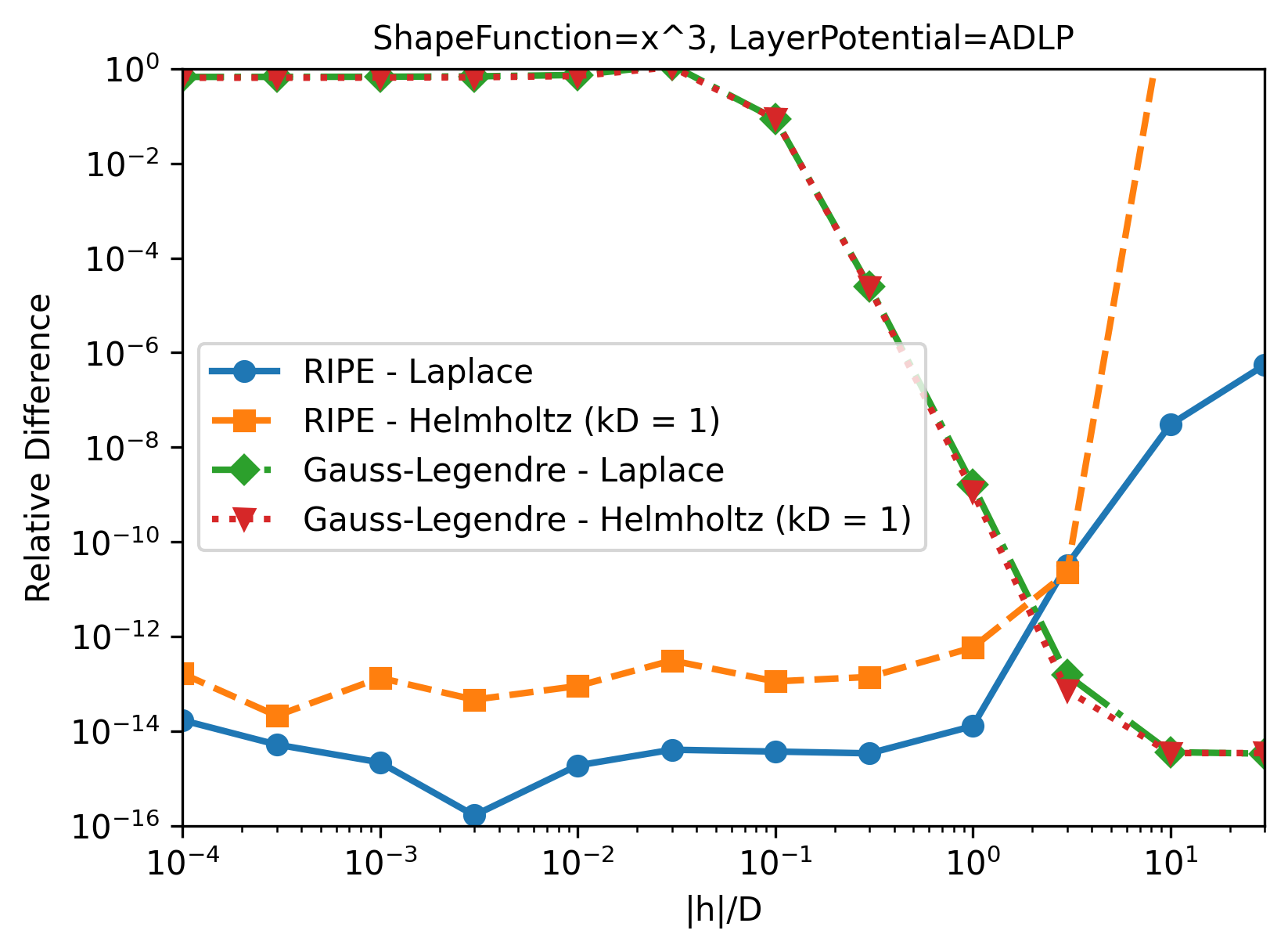}
    \includegraphics[width=6cm]{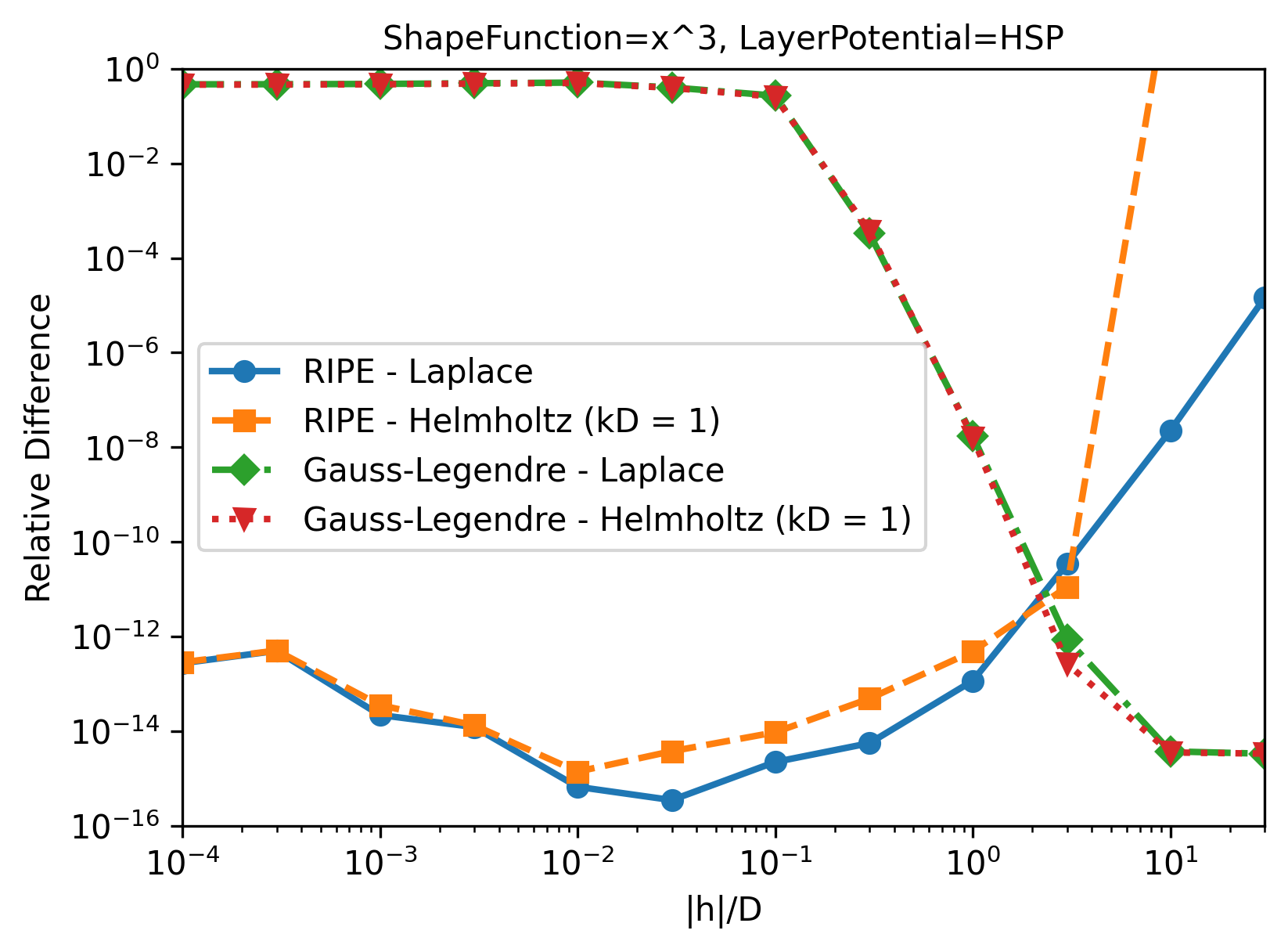}
    \caption{Relative difference of the layer potentials between RIPE or 12-th order Gauss-Legendre quadrature and the reference adaptive Gauss-Kronrod quadrature, for the shape function $N(x,y)=x^3$.
            Results for single layer potential (top left), double layer potential (top right), adjoint double layer potential (bottom left), and hypersingular potential (bottom right).
            }
    \label{fig:recursiveIntegrals_xto3}
    \end{figure*}

The method was tested for all four layer potentials, for both the Laplace and Helmholtz kernels. Adaptive Gauss-Kronrod quadrature, implemented in QUADPACK~\cite{piessens2012quadpack}, was used to compute the  layer potentials $P_{\mathrm{GK}}$ over a boundary element. The result of RIPE, $P_{\mathrm{RIPE}}$ was compared against this reference result via the relative error $|P_{\mathrm{RIPE}}-P_{\mathrm{GK}}|/|P_{\mathrm{GK}}|$. The maximum expansion order $p_e$ in RIPE was set to 32, and the error tolerance of Gauss-Kronrod was set to $10^{-12}$. 
The layer potentials were also evaluated using 12th order Gauss-Legendre quadrature which has 33 quadrature nodes on a triangle element. 
The triangle $\{(0,0,0),(1,0,0),(0,1,0)\}$ and point $\rp=(1/3,1/3,|h|)$ were used as the element and observation points, respectively. 
$N(x,y)=x^3$ was used as the shape function. 
Results are shown in \cref{fig:recursiveIntegrals_xto3}. 
Our RIPE method delivers remarkable accuracy particularly in the challenging nearly singular regime $|h|/D<1$, with $D$ the maximum edge length of the element.  
Gauss-Legendre quadrature is not reliable in this domain. 
Note that the condition $kD=1$ approximately corresponds to six wavelengths per element, which is typically used as the maximum mesh size in boundary element analysis. 
We also observe that the RIPE method has worse performance for larger  $|h|/D$, where Gauss-Legendre, on the other hand, delivers good accuracy. 
Given this complementary character, one can switch the integration routine based  on the distance of the evaluation point from the element. 

\subsubsection{Singular and hypersingular case}
Cases with the observation point $\rp$ on the element (singular and hypersingular cases) are handled without any modifications to the framework. 
For these the method was compared with results obtained by Guiggiani's method~\cite{guiggiani1992general}, internally using Gauss-Legendre quadrature of 20th order. 
The same element as in the previous section with $\rp$ on the center of the element was used to compute the single layer and the hypersingular potentials. The six shape functions of a conforming second-order Lagrange triangle element were used. 
Results for both Laplace and Helmholtz kernels are shown in \cref{table_singular_laplace,table_hypersingular_laplace,table_singular_helmholtz,table_hypersingular_helmholtz}. 
A relative difference $|P_{\mathrm{RIPE}} - P_{\mathrm{Gui}}|/|P_{\mathrm{Gui}}|$ from $10^{-10}$ to $10^{-13}$ was observed. 
\begin{table*}[htbp]
\centering
\begin{tabular}{|c|c|c|c|}
\toprule\hline
{} &  $P_{\mathrm{RIPE}}$ &  $P_{\mathrm{Gui}}$ &  $|P_{\mathrm{RIPE}}-P_{\mathrm{Gui}}|/|P_{\mathrm{Gui}}|$\\
\midrule\hline
$N_1$ &   -0.0096108650741614 &       -0.0096108650753968 &  $1.28540284\times 10^{-10}$ \\
$N_2$ &   -0.0096108650741614 &       -0.0096108650753968 &  $1.28539201\times 10^{-10}$ \\
$N_3$ &   -0.0059161308348599 &       -0.0059161308348537 &  $1.05940117\times 10^{-12}$ \\
$N_4$ &    0.0733163156462961 &        0.0733163156487268 &  $3.31537174\times 10^{-11}$ \\
$N_5$ &    0.0716914080260122 &        0.0716914080259944 &  $2.48939665\times 10^{-13}$ \\
$N_6$ &    0.0716914080260122 &        0.0716914080259944 &  $2.48552512\times 10^{-13}$ \\
\bottomrule\hline
\end{tabular}
\caption{Single layer potentials for the Laplace kernel evaluated over the element with second-order shape functions of a conforming Lagrange triangle element. Numbers for RIPE ($P_{\mathrm{RIPE}}$), Guiggiani's method ($P_{\mathrm{Gui}}$), and their relative differences are shown.}
\label{table_singular_laplace}
\end{table*}
\begin{table*}[htbp]
\centering
\begin{tabular}{|c|c|c|c|}
\toprule\hline
{} &  $P_{\mathrm{RIPE}}$ &  $P_{\mathrm{Gui}}$ &  $|P_{\mathrm{RIPE}}-P_{\mathrm{Gui}}|/|P_{\mathrm{Gui}}|$ \\
\midrule\hline
$N_1$ &    0.3411586129005689 &        0.3411586129009882 &  $1.22913865\times 10^{-12}$ \\
$N_2$ &    0.3411586129005690 &        0.3411586129009879 &  $1.22799965\times 10^{-12}$ \\
$N_3$ &    0.5031187119584526 &        0.5031187119589535 &  $9.95654936\times 10^{-13}$ \\
$N_4$ &   -0.7261344637586460 &       -0.7261344637604444 &  $2.47674552\times 10^{-12}$ \\
$N_5$ &   -0.9322819538428125 &       -0.9322819538448199 &  $2.15320509\times 10^{-12}$ \\
$N_6$ &   -0.9322819538428122 &       -0.9322819538448149 &  $2.14820345\times 10^{-12}$ \\
\bottomrule\hline
\end{tabular}
\caption{Hypersingular potentials for the Laplace kernel evaluated over the second-order element.}
\label{table_hypersingular_laplace}
\end{table*}
\begin{table*}[htbp]
\centering
\begin{tabular}{|c|c|c|c|}
\toprule\hline
{} &  $P_{\mathrm{RIPE}}$ &  $P_{\mathrm{Gui}}$ &  $|P_{\mathrm{RIPE}}-P_{\mathrm{Gui}}|/|P_{\mathrm{Gui}}|$ \\
\midrule\hline
$N_1$ &   -0.0097575874677327 &       -0.0097575874673111 &  $1.65408707\times 10^{-10}$ \\
$N_2$ &   -0.0097575874677327 &       -0.0097575874673111 &  $1.65411332\times 10^{-10}$ \\
$N_3$ &   -0.0059358291069226 &       -0.0059358291069324 &  $2.15796857\times 10^{-12}$ \\
$N_4$ &    0.0724350497721009 &        0.0724350497713923 &  $4.33315252\times 10^{-11}$ \\
$N_5$ &    0.0707990955161933 &        0.0707990955162131 &  $4.29383222\times 10^{-13}$ \\
$N_6$ &    0.0707990955161934 &        0.0707990955162131 &  $4.30321106\times 10^{-13}$ \\
\bottomrule\hline
\end{tabular}
\caption{Single layer potentials for the Helmholtz kernel ($kD=1$) over the second-order element. }
\label{table_singular_helmholtz}
\end{table*}
\begin{table*}[htbp]
\centering
\begin{tabular}{|c|c|c|c|}
\toprule\hline
{} &  $P_{\mathrm{RIPE}}$ &  $P_{\mathrm{Gui}}$ &  $|P_{\mathrm{RIPE}}-P_{\mathrm{Gui}}|/|P_{\mathrm{Gui}}|$ \\
\midrule\hline
$N_1$ &    0.3387374371700406 &        0.3387374371703727 &  $1.34783671\times 10^{-12}$ \\
$N_2$ &    0.3387374371700407 &        0.3387374371703721 &  $1.34571877\times 10^{-12}$ \\
$N_3$ &    0.5016372264001558 &        0.5016372264006584 &  $1.00193060\times 10^{-12}$ \\
$N_4$ &   -0.7079157406214216 &       -0.7079157406230446 &  $2.45050302\times 10^{-12}$ \\
$N_5$ &   -0.9144708322784949 &       -0.9144708322805092 &  $2.20267780\times 10^{-12}$ \\
$N_6$ &   -0.9144708322784947 &       -0.9144708322805069 &  $2.20036963\times 10^{-12}$ \\
\bottomrule\hline
\end{tabular}
\caption{Hypersingular potentials for the Helmholtz kernel ($kD=1$) over the second-order element. }
\label{table_hypersingular_helmholtz}
\end{table*}
\subsection{Stability}
\begin{figure*}[htb!]
    \centering
    \includegraphics[width=6cm]{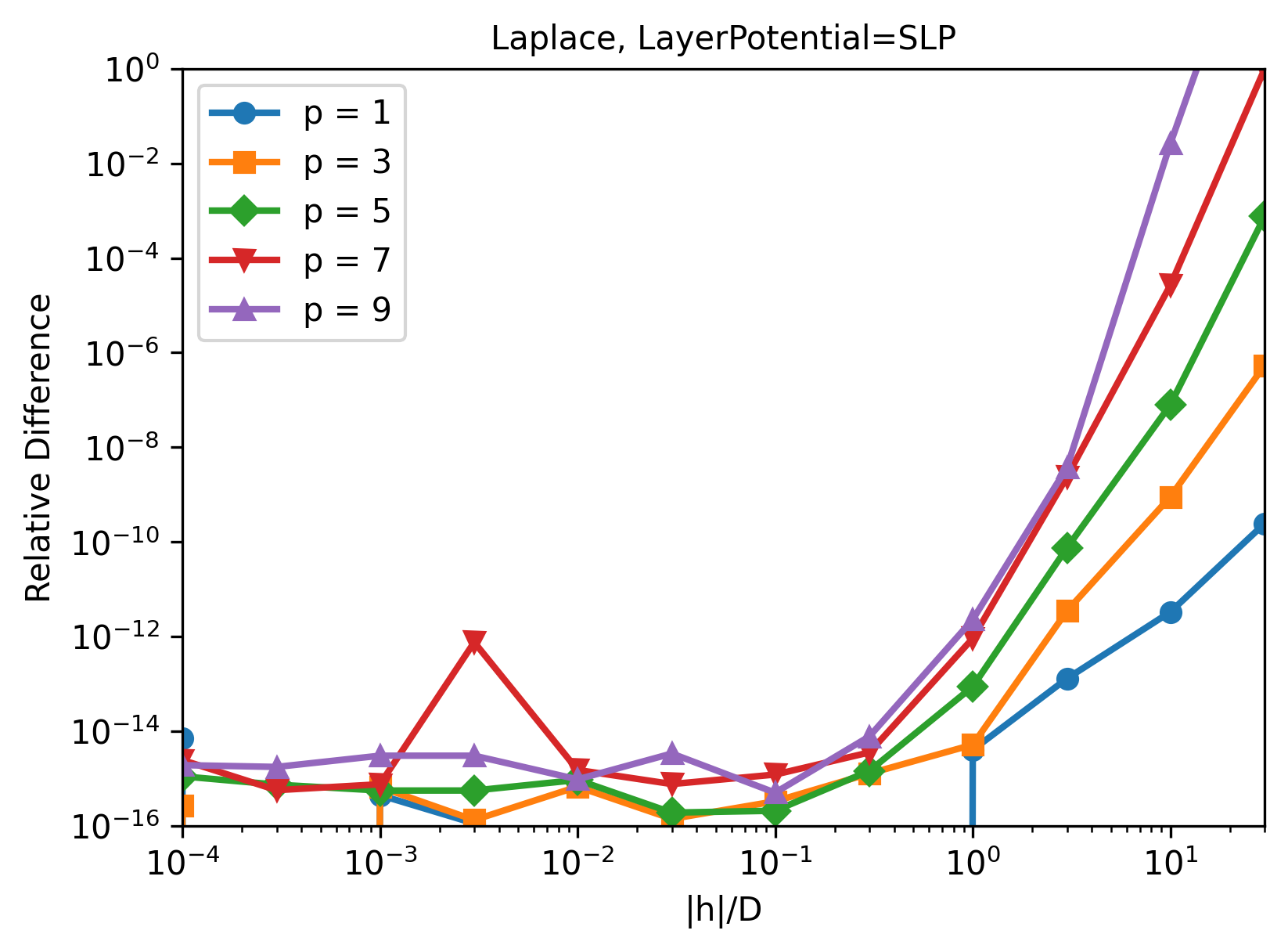}
    \includegraphics[width=6cm]{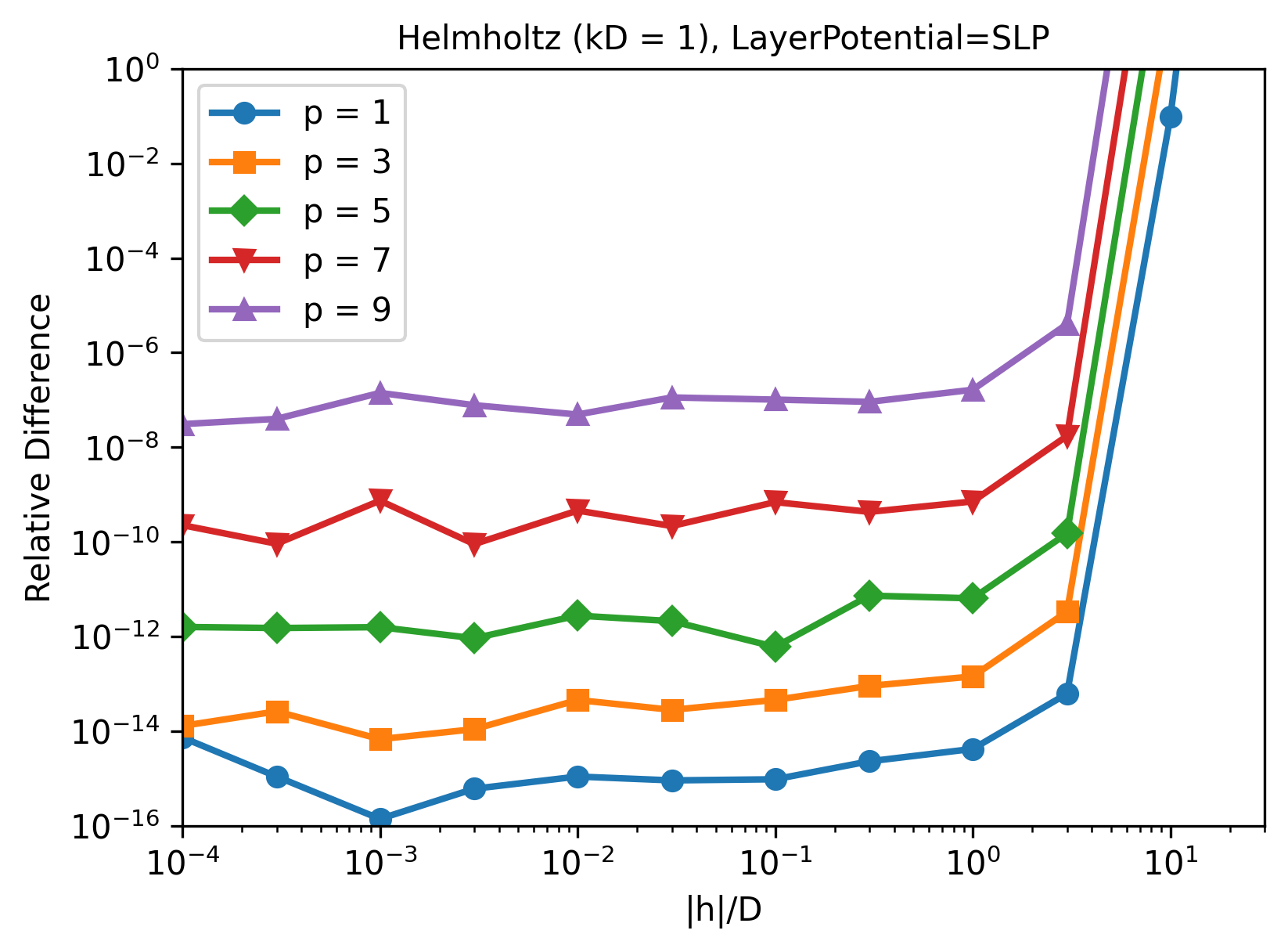}
    \includegraphics[width=6cm]{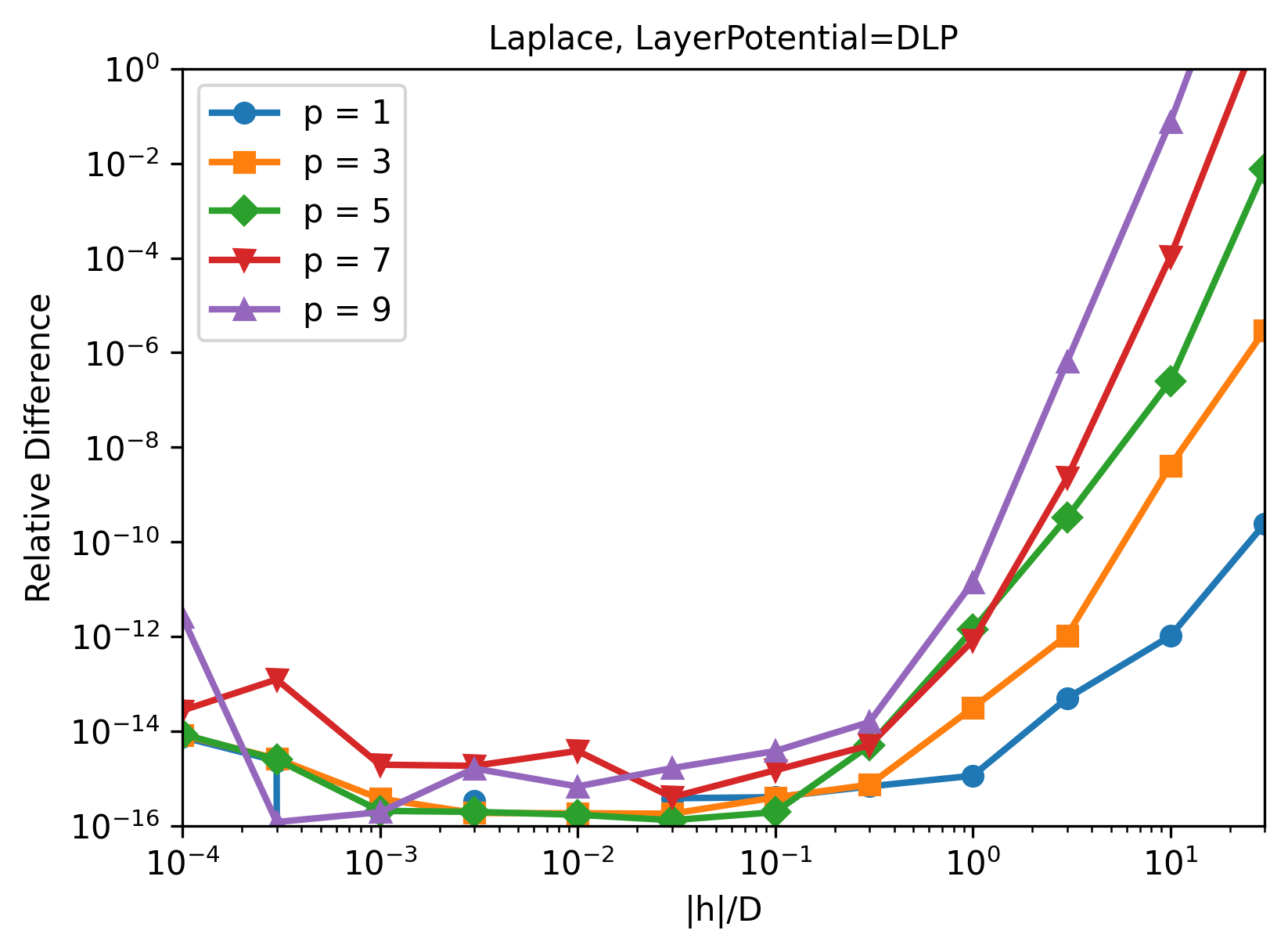}
    \includegraphics[width=6cm]{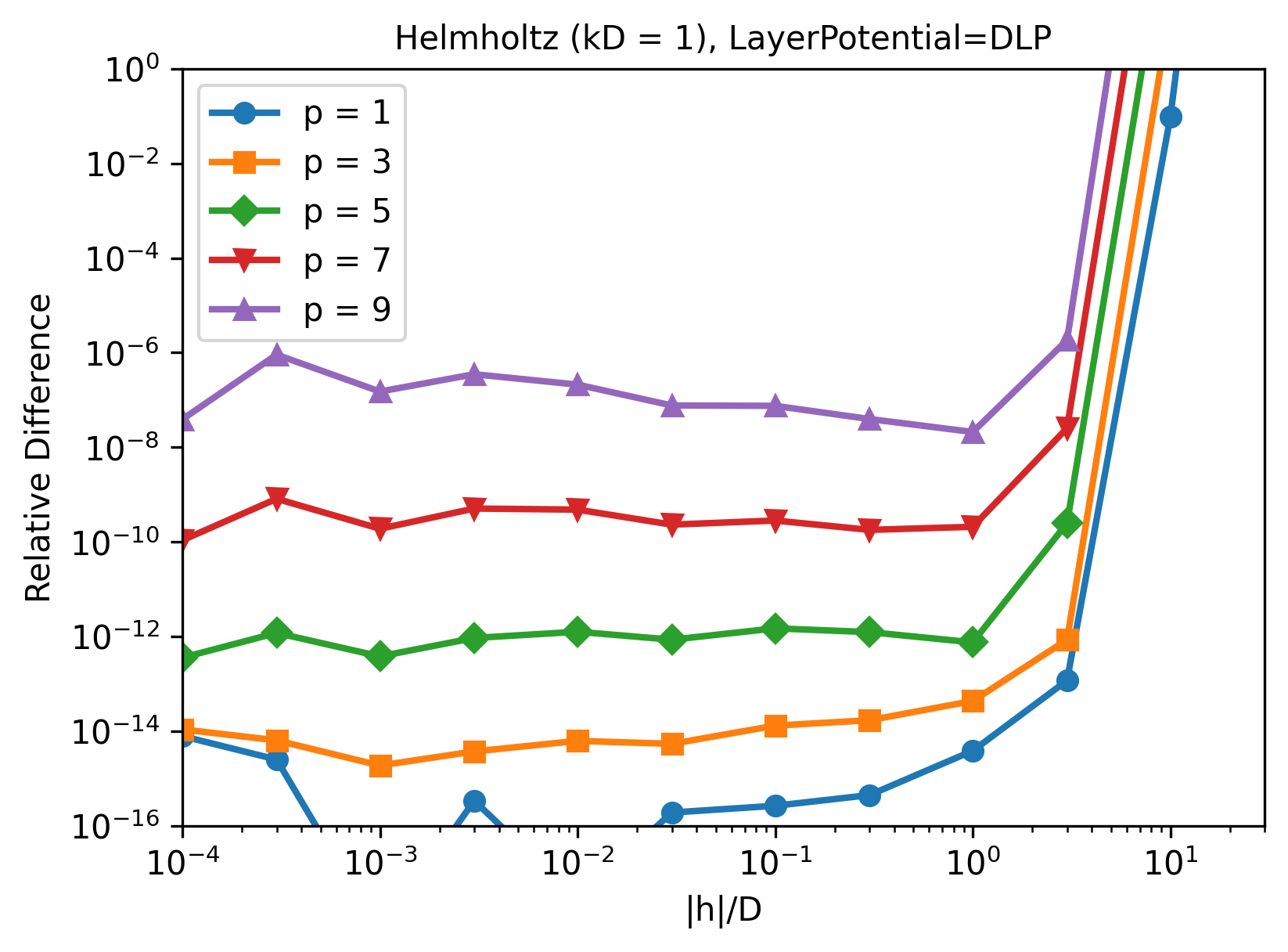}
    \includegraphics[width=6cm]{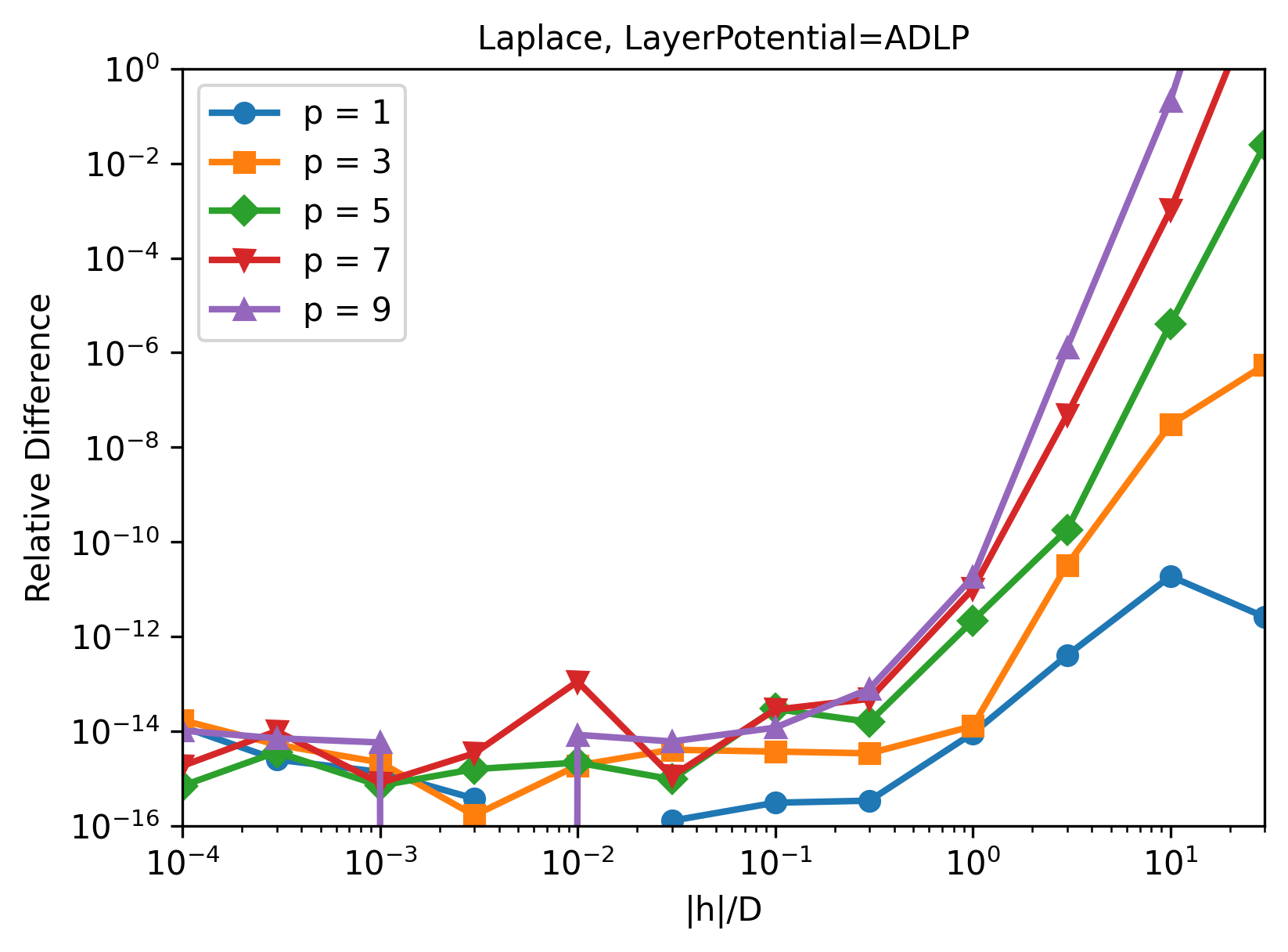}
    \includegraphics[width=6cm]{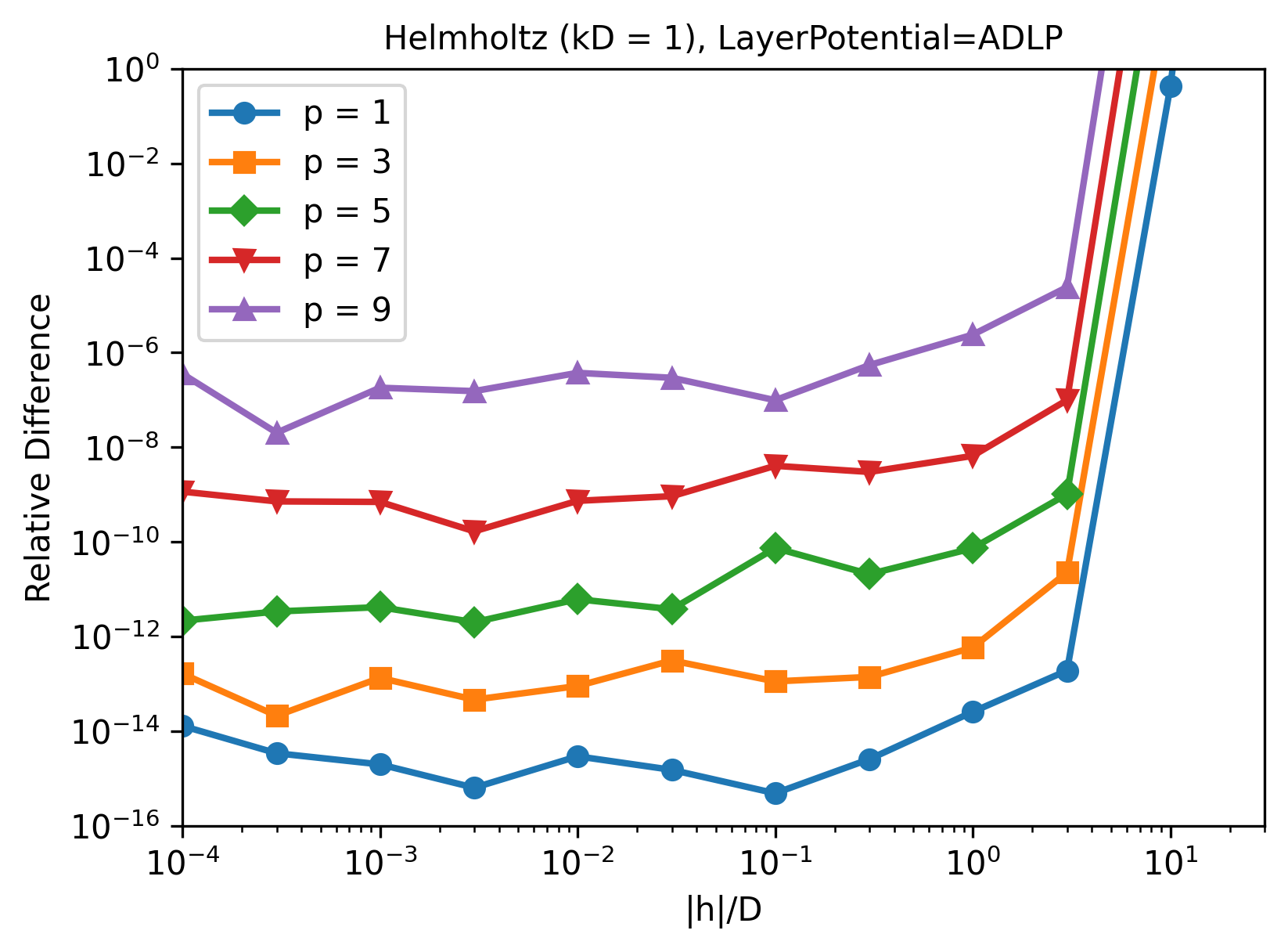}
    \includegraphics[width=6cm]{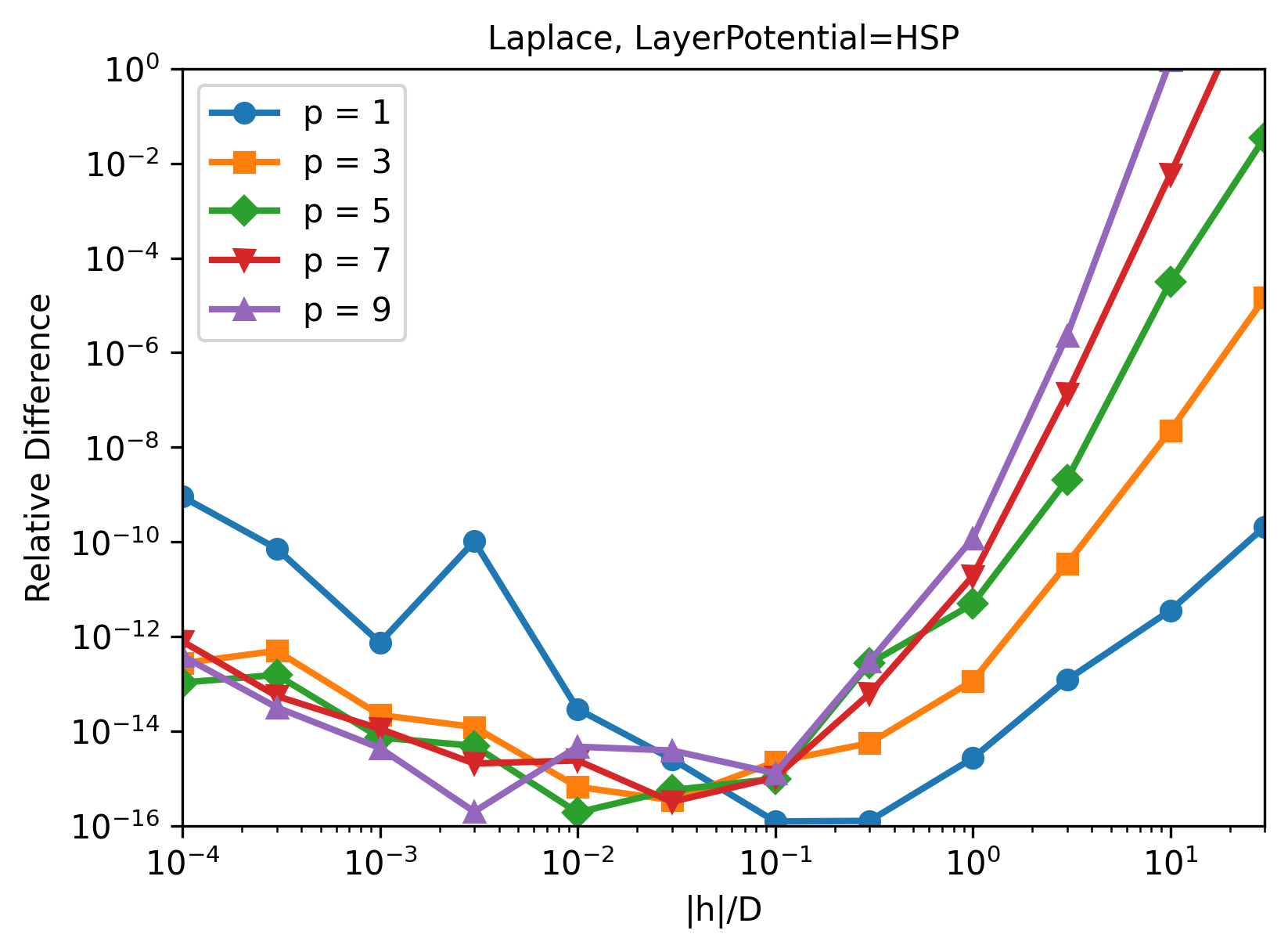}
    \includegraphics[width=6cm]{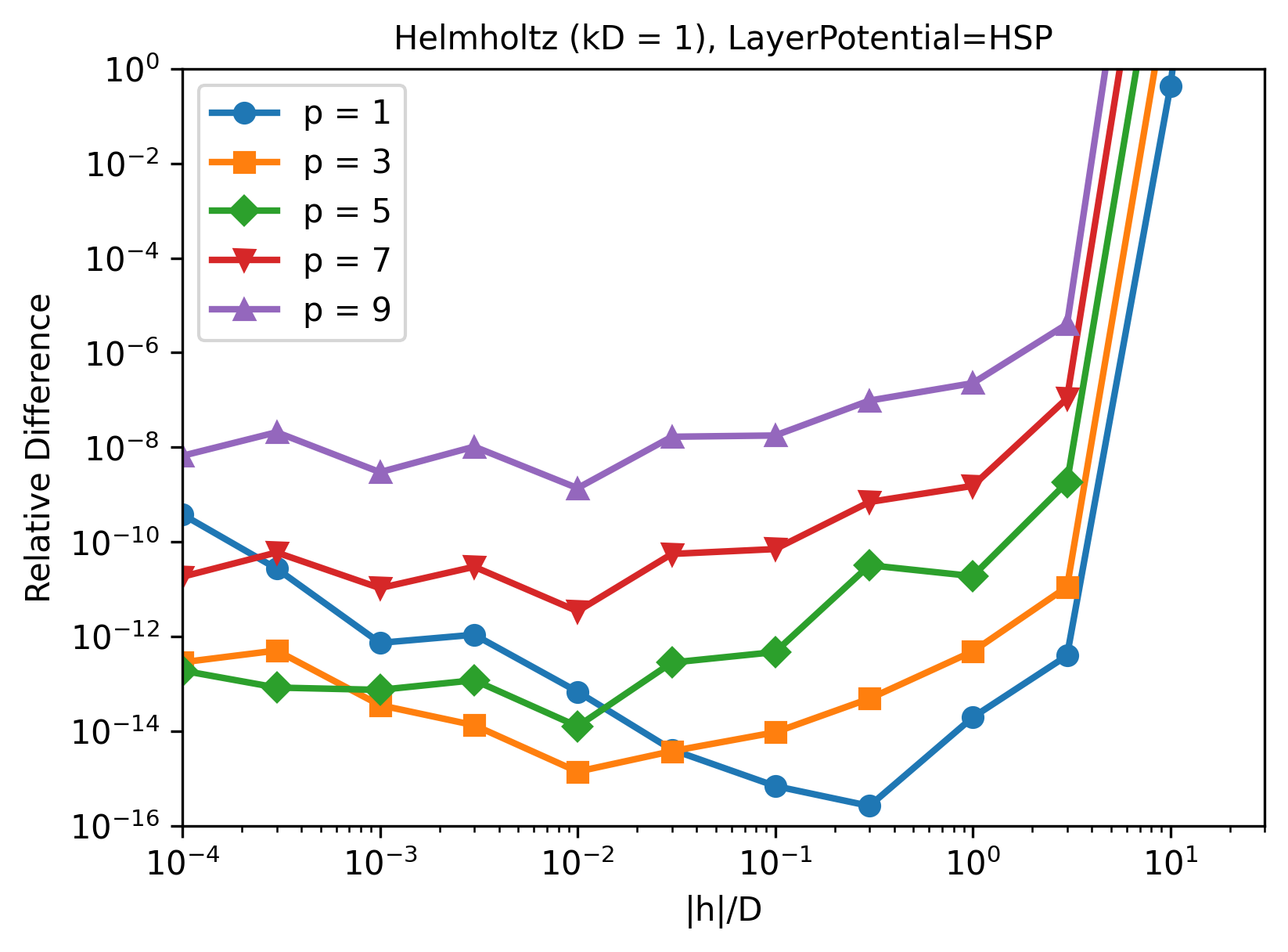}
    \caption{Relative difference of the layer potentials between RIPE and the adaptive Gauss-Kronrod quadrature, for various shape function orders. 
    Laplace (left column) and Helmholtz (right column), for single layer, double layer, adjoint double layer, and hypersingular potentials (top to bottom).}
    \label{fig:recursiveIntegrals_stability}
    \end{figure*}

The RIPE method was evaluated by testing its accuracy with high order elements up to $p_s=9$. 
\cref{fig:recursiveIntegrals_stability} shows results for both Laplace and Helmholtz layer potentials. 
The relative difference with the reference Gauss-Kronrod increases with element order, and this divergence is more for the Helmholtz kernel.
\subsection{Computation time}
\begin{figure}[htbp]
\centering
\includegraphics[width=8.5cm]{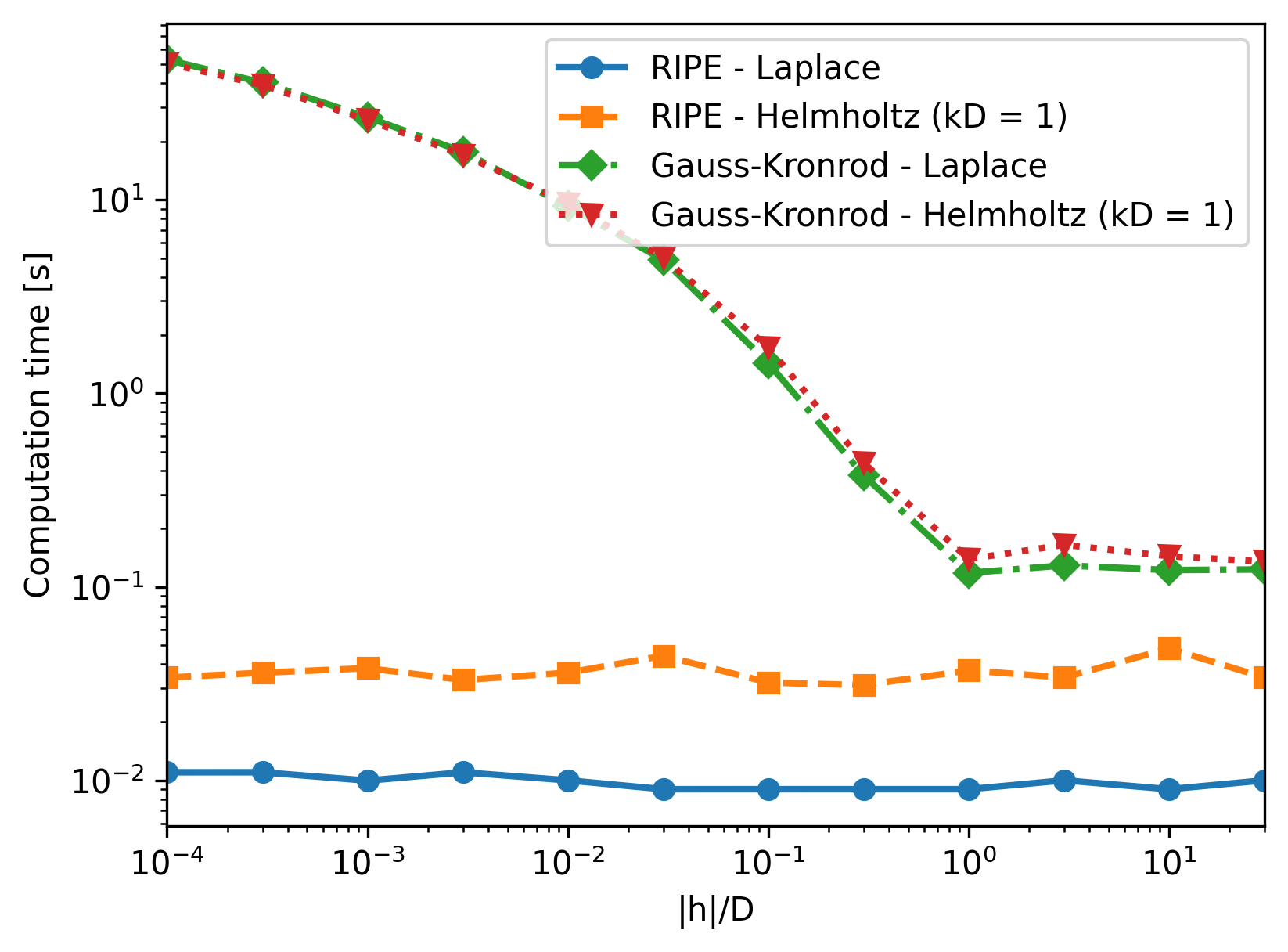}
\caption{The total wall-clock time for computing all four layer potentials for a 5th order element, for RIPE and for adaptive Gauss-Kronrod quadrature. Note that the times for RIPE are the times spent to compute integrals weighted by $x^b y^c$ for all $((p_s+1)(p_s+2))/2=21$ terms in $0\le b+c \le 5$ in one run of the recursion, while the times for Gauss-Kronrod are for computing only the integral with shape function $x^{5}$. 
}
\label{fig:times}
\end{figure}

Computation time of the layer potentials with 5th order shape functions over the element $\{(0,0,0),(1,0,0),(0,1,0)\}$ and observation point $\rp=(1/3,1/3,|h|)$ were measured and compared against the adaptive Gauss-Kronrod baseline. 
The expansion order for the Helmholtz kernel in RIPE was set to $p_e=32$, and the error tolerance of Gauss-Kronrod was set to $10^{-12}$.
The result is shown in \cref{fig:times}.
The times for RIPE are those spent to compute integrals weighted by the shape function $N(x,y)=x^b y^c$ for all terms in $0\le b+c \le p_s$ in one run of the recursion, while for Gauss-Kronrod these are for the computation of a {\em single} term with the shape function $N(x,y)=x^{5}$. 
Gauss-Kronrod requires long computation time in the near-singular domain with small $|h|/D$, while computation time for RIPE does not depend on $|h|/D$ and can be orders of magnitudes faster. 
The measurements for RIPE were performed using a Python prototype without performance optimization, and Gauss-Kronrod was measured using Scipy's interface to QUADPACK. 
The advantage of the RIPE method should be more significant in a production collocation BEM setup where the $x^b y^c$ terms for all the $0 \le b+c \le p_s$ combinations are needed.

\section{Conclusion}
A recursive algorithm (RIPE) to evaluate layer potentials arising in the collocation BEM for the Laplace and Helmholtz equation, tailored for piecewise flat boundary elements with polynomial shape functions of arbitrary orders was proposed. 
Numerical tests showed that RIPE exhibits remarkable accuracy in the nearly-singular regime, where Gauss-Legendre quadrature is not effective. 
Experiments indicate that RIPE is several orders of magnitude faster than Gauss-Kronrod quadrature, a general-purpose adaptive quadrature method. 
RIPE would hence serve as an efficient routine for nearly singular, singular, and hypersingular integrals for Laplace and Helmholtz layer potentials. 
The benefits of the RIPE method are that it offers: (1) analytical integration for polynomial elements achieved by a simple formulation using auxiliary vector fields, (2) easier error control for the Helmholtz case compared to methods based on Gauss-Legendre quadrature due to the series expansion-based formulation, (3) accuracy and efficiency for nearly singular, singular, and hypersingular integrals, achieved under a single framework naturally supporting all of these cases without separate modifications, for all four standard layer potentials, for both Laplace and Helmholtz kernels. 
 
The RIPE method has multiple frontiers for further development. 
While it is restricted to flat elements, supporting manifold surfaces as shown in \cite{zhu2022high} is important direction for development since the reduction of geometrical errors via higher-order geometry representations can be beneficial for curved objects~\cite{greengard2021fast}.
We note, however, that many practical BEM software are written for flat elements, and many problems have meshes from architecture or CAD that are truly flat, or where plane triangular meshes are more easily available. 
Extending RIPE to Galerkin BEM, other kernels, and integrating it into FMM-BEM solvers are other directions for future work.

\appendix
\section{Elementary integral computations}\label{appendix_kl_il}
\subsection{Computation of integrals $i_m$}
We show analytical computation of 
\begin{equation}
i_{m}\left( x;a\right) =\int r^{m}dx,\quad r=\sqrt{x^{2}+a^{2}},\quad
a^{2}=y^{\prime 2}+z^{\prime 2},\quad m=0,\pm 1,\pm 2,...  \label{A1}
\end{equation}
For small values of even and odd $m$ these
integrals can be computed 
\begin{eqnarray}
i_{0}\left( x;a\right) &=& x,\quad i_{-1}\left( x;a\right) = \ln \left| r+x\right| ,\quad  \label{A1.1}
i_{-2}\left( x;a\right) = \frac{1}{a}\arctan \frac{x}{a}%
,\quad i_{-3}\left( x;a\right) = \frac{x}{a^{2}r}.  \nonumber
\end{eqnarray}
We have the recurrence: 
\begin{eqnarray}
i_{m+2} &=&\int r^{m+2}dx=xr^{m+2}-\left( m+2\right) \int x^{2}r^{m}dx
\label{A1.2} 
= xr^{m+2}-\left( m+2\right) i_{m+2}+(m+2)a^{2}i_{m}.  \nonumber
\end{eqnarray}
\begin{equation}
i_{m+2}=\frac{xr^{m+2}}{m+3}+\frac{m+2}{m+3}a^{2}i_{m},\quad m\neq -3.
\label{A1.3}
\end{equation}
We need integral values for $m\geqslant -3$, and have explicit expressions for all required non-positive $m $ and can recursively find all positive $m$ starting the recurrence from $m=0$ for even $m$ and from $%
m=-1$ for odd $m$.  
Note that primitives $i_{m}\left( x;a\right) $ may have singularities if $a=0.$ For these we have 
\begin{equation}
i_{m}\left( x;0\right) =\frac{x\left| x\right| ^{m}}{m+1},\quad m\neq
-1,\quad i_{-1}\left( x;0\right) =\mbox{sgn}(x)\ln \left| x\right| .  \label{A6.1}
\end{equation}%
\subsection{Computation of integrals $k_m$}
In this appendix we show how integrals 
\begin{equation}
k_{m}\left( x;y^{\prime },z^{\prime }\right) =z^{\prime }\int \frac{r^{m}}{x^{2}+z^{\prime 2}}dx,\quad r=\sqrt{x^{2}+y^{\prime 2}+z^{\prime 2}},\quad
m=0,\pm 1,\pm 2,...  \label{B1}
\end{equation}%
can be computed analytically. Note, that for $y^{\prime }=0$ these reduce to
the integrals $i_{m}$: 
\begin{equation}
k_{m}\left( x;0,z^{\prime }\right) =z^{\prime }\int \frac{r^{m}}{r^{2}}%
dx=z^{\prime }i_{m-2}\left( x;\left| z\right| \right) .  \label{B1.0}
\end{equation}%
So we only need consider the case $y^{\prime }\neq 0.$ We derive recurrences
for all needed $k_{m}$ and find initial values to use them. We have 
\begin{eqnarray}
k_{m+2} &=&z^{\prime }\int \frac{r^{m+2}}{x^{2}+z^{\prime 2}}dx=z^{\prime
}\int \frac{\left( x^{2}+z^{\prime 2}\right) +y^{\prime 2}}{x^{2}+z^{\prime
2}}\rho ^{m}dx  \label{B1.1} 
= z^{\prime }i_{m}\left( x;\sqrt{y^{\prime 2}+z^{\prime 2}}\right)
+y^{\prime 2}k_{m}.  \nonumber
\end{eqnarray}%
With known expressions for the initial values (odd and even), all integrals
can be computed recursively. For this we have 
\begin{eqnarray}
k_{-1}\left( x;y^{\prime },z^{\prime }\right) &=& \frac{z^{\prime }}{y^{\prime }\left| z^{\prime }\right| }\arctan \frac{%
y^{\prime }x}{\left| z^{\prime }\right| r}, \quad
k_{0}\left( x;y^{\prime },z^{\prime }\right)  = \frac{z^{\prime }}{\left|
z^{\prime }\right| }\arctan \frac{x}{\left| z^{\prime }\right| }, \\
k_{1}\left( x;y^{\prime },z^{\prime }\right)  &=& \frac{y^{\prime }z^{\prime }%
}{\left| z^{\prime }\right| }\arctan \frac{y^{\prime }x}{\left| z^{\prime
}\right| r}+z^{\prime }\ln \left| r+x\right|
. \nonumber
\end{eqnarray}%
Note that $y^{\prime }k_{-1}\left( x;y^{\prime },z^{\prime }\right) $
entering the primitive expressions are not singular even when $z^{\prime }$
and $y^{\prime }$ approaches zero, while depend on the path (ratio $%
y^{\prime }/z^{\prime }$).

\bibliographystyle{siamplain}

\end{document}

%% file: ex_shared.tex
\usepackage{lipsum}
\usepackage{amsfonts,amssymb,booktabs}
\usepackage{graphicx}
\usepackage{epstopdf}
\usepackage{algorithmic}
\ifpdf
  \DeclareGraphicsExtensions{.eps,.pdf,.png,.jpg}
\else
  \DeclareGraphicsExtensions{.eps}
\fi

\newsiamremark{remark}{Remark}
\newsiamremark{hypothesis}{Hypothesis}
\crefname{hypothesis}{Hypothesis}{Hypotheses}
\newsiamthm{claim}{Claim}

\headers{Recursive Quadrature of Layer Potentials}{S. Kaneko, N. A. Gumerov, and R. Duraiswami}

\title{Recursive Analytical Quadrature of Laplace and Helmholtz Layer Potentials 
in $\mathbb{R}^3$ \thanks{Submitted to the editors December 6th, 2022.
\funding{This work is partially supported by Cooperative Research Agreement  W911NF2020213 between the University of Maryland and the Army Research Laboratory, with David Hull, Ross Adelman
and Steven Vinci as Technical monitors. 
Shoken Kaneko received scholarships from Japan Student Services Organization and Watanabe Foundation.
}}}

\author{Shoken Kaneko\thanks{Department of Computer Science, University of Maryland, College Park, MD (\email{kaneko60@umd.edu}). }
\and Nail A. Gumerov\thanks{UMIACS, University of Maryland, College Park, MD (\email{ramanid@umd.edu}).}
\and Ramani Duraiswami $^\dagger$ \footnotemark[3]
}

\usepackage{amsopn}
